\documentclass[a4paper,10pt,english]{article}

\usepackage[T1]{fontenc}
\usepackage{times}

\usepackage{hyperref}

\usepackage[centertags]{amsmath}
\usepackage{amsthm}
\usepackage{amssymb}
\usepackage{bbold} % Questo pacchetto non viene trovato anche se e' effettivamente installato
\usepackage{verbatim}

\newcounter{thmcount}
\numberwithin{thmcount}{section}
\theoremstyle{plain}
\newtheorem{Theorem}[thmcount]{Theorem}
\newtheorem{Proposition}[thmcount]{Proposition}
\newtheorem{Lemma}[thmcount]{Lemma}
\newtheorem{Corollary}[thmcount]{Corollary}

\theoremstyle{definition}
\newtheorem{Definition}[thmcount]{Definition}

\theoremstyle{remark}
\newtheorem{Remark}[thmcount]{Remark}
\newtheorem{Example}[thmcount]{Example}

\usepackage{graphicx}
\usepackage{varioref}

\newcommand\defeq{{\stackrel{\mbox{\tiny{def}}}{=}}}

\def\ammark#1{}%\texttt{[AM:#1]}}
\def\admark#1{}%\texttt{[AD:#1]}}

\makeatletter
\def\len{\mathop{\operator@font len}\nolimits}
\def\dil{\mathop{\operator@font dil}\nolimits}
\def\argmin{\mathop{\operator@font argmin}\nolimits}
\makeatother

\def\e{\varepsilon}

\newcommand{\uno}{{\mathbb 1}}
\newcommand\interno[1]{{\mathaccent"7017 {#1}}}

\newcommand{\calI}{\mathcal{I}}

\newcommand{\real}{{\mathbb R}} % Numeri Reali
\newcommand{\N}{\mathbb N}    % Numeri Naturali Positivi

\newcommand{\qbrk}[1]{\left[ {#1} \right]}   % parentesi quadre
\newcommand{\gbrk}[1]{\left\{ {#1} \right\}} % parentesi graffe
\newcommand{\lbrk}[1]{\left| {#1} \right|} 	 % parentesi linea
\newcommand{\set}[1]{\gbrk{#1}}

\begin{document}

\title{Banach-like metrics and metrics of compact sets.}
\author{A. Duci, A. C. Mennucci}
\date{\today}
\maketitle

%%%%%%%%%%%%%%%%%%%%%%%%%%%%%%%%%%%%%%%%%%%%%%%%%%%%%%%%%%%%%%%%%%%%%%%%%%%%%%%%%%%%%%%%%% 

\begin{abstract}
  We present and study a family of metrics on the space of compact
  subsets of $\real^N$ (that we call ``shapes'').  These metrics
  are ``geometric'', that is, they are independent of
  rotation and translation; and these metrics enjoy many interesting
  properties, as, for example, the existence of minimal geodesics.
  We view our space of shapes as a subset of Banach (or Hilbert) manifolds:
  so we can   define a ``tangent manifold'' to shapes,
  and (in a very weak form) talk of a  ``Riemannian Geometry''  of shapes.
  Some of the metrics that we propose are topologically equivalent to the
  Hausdorff metric; but at the same time, they are more ``regular'', since 
  we can hope for a local uniqueness of minimal geodesics.

  We also study properties of the metrics obtained by isometrically
  identifying a generic metric space with a subset of a Banach space to obtain a
  rigidity result.
\end{abstract}


\noindent
{\bf Acknowledgments.} A. Duci was supported by the EU 6th Framework Programme Grant MRTN - CT - 2005 - 019481.

\section{Introduction}

A wide interest for the study of \emph{shape spaces} arose in recent
years, in particular inside the Computer Vision community.
%and the image processing community. 

There are two different (but interconnected) fields of applications for
a good Shape Space in Computer Vision: 
\begin{description}
\item[Shape Optimization] where we want to find the shape that best 
  satisfies a design goal;
  a topic  interest in Engineering at large;
\item[Shape Analysis] where we study a family of Shapes
  for purposes of  statistics, (automatic) cataloging, probabilistic modeling,
  among others, and possibly 
  to create an a-priori model for a better Shape Optimization.
\end{description}

To achieve the above, some structure is clearly needed on the Shape
Space, so that our goals can be studied and the problem can be solved.

\begin{Remark}\label{rem:preshape}
\em  Note that, for the purpose of Shape Optimization, shapes
  are usually intended ``up to    rotation, translation and scaling'';
  for this reason, when we wish to distinguish between the two,  we will call a space
  for Shape Optimization a ``\emph{preshape space}''.
\end{Remark}

\subsection{Shape spaces}
In general the ``Shape Space'' $\calI$ will be  a suitable choice  of subsets of $\real^N$.

A common way to model shapes is by 
\textbf{representation/embedding}:
\begin{itemize}
\item we  \textbf{represent} the shape $A$ by a function  $u_A$
\item and then we \textbf{embed} this representation in
  a space $E$, so that we can operate on the 
  shapes $A$ by operating on the representations $u_A$;
\end{itemize}
for example, if $E$ is a Banach space with norm $\|\cdot\|$,
we can define a \emph{distance of shapes}
by $d(A,B)\defeq \| u_A - u_B\|$.

Most often, this representation/embedding scheme does not 
directly provide a Shape Space satisfying all desired properties.
In particular, in many cases it happens that the representation
is ``redundant'', that is, the same shape has many different possible
representations. An appropriate \textbf{quotient} is then introduced.

% Two possible example definitions of ``Shape Space'' $\calI$ are 
% \begin{itemize}
% \item the space of  (embedded) curves; or, more in general,
% \item a suitable choice  of subsets of $\real^N$;
% \end{itemize}
% (and, in this paper we will indeed deal with the space of all compact sets).

\smallskip

There are many examples of the \emph{representation/embedding/quotient} scheme in the literature;
for the case of generic subsets of $\real^N$,
\begin{itemize}
\item a standard  representation is obtained 
  by associating a closed subset $A$ to the \textbf{distance function}
  \begin{equation}
    u_A(x)\defeq \inf_{y\in A} |x-y| \label{eq:distfun}
  \end{equation}
  or the \textbf{signed distance function}
  \begin{equation}
    b_A(x)\defeq u_A(x) - u_{\real^N\setminus A}(x)\label{eq:signdistfun}
  \end{equation}

  We can then define a \emph{topology of shapes} by deciding that
  $A_n\to A$ when $u_{A_n}\to u_A$  uniformly on compact sets. This convergence coincides with
  the Kuratowski topology of closed sets.

  We can also operate \emph{``linearly''} on shapes by operating on $u_A$ or $b_A$:
  so  we can define \emph{shape averages} and \emph{shape principal component analysis}.
  Note that in general a linear combination of (signed) distance functions will not
  be a (signed) distance function:
  so any linear operation must be followed by an \emph{ad hoc} correction.
  For example, given two shapes $A_0,A_1$, we can define an interpolation  $A_t$ for $t\in[0,1]$
  by setting  $A_t = \{ x ~|~ t b_{A_1}(x)+ (1-t) b_{A_0}(x)\le 0\}$.

  This Shape Space is not independent of the position: when it is used
  for shape analysis, a \emph{registration} of the shapes to a common
  pose is often performed (but, see also sec. \ref{sec:dist_quot}).

\item A. Duci \emph{et al} (see \cite{duci03:harmon_embed,
    duci06:harmon_embed}) represent a closed planar contour as the
  zero level of a harmonic function. This novel representation for
  contours is explicitly designed to possess a linear structure, which
  greatly simplifies linear operations such as averaging, principal
  component analysis or differentiation in the space of shapes.

\item Trouv\'e--Younes \emph{et al} (see \cite{GTY05:modeling},
  \cite{MR2176922} and references therein) modeled the motion of
   shapes by studying a left invariant Riemannian metric on the
  diffeomorphisms of the space $\real^n$; to recover a true metric of
  shapes, a quotient is then added.
\end{itemize}

But the \emph{representation/embedding/quotient} scheme is also found
when dealing with spaces of curves:
\begin{itemize}
\item In the work of Mio, Srivastava \emph{et al}.
  \cite{Srivastava:StatAna,mio_sriv_04:elast_strin_model,Srivastava:AnalPlan},
  smooth planar closed curves
  $c:S^1\to\real^2 $ of length $2\pi$ are  parametrized by arclength
  and represented by the angle function $\alpha
  [0,2\pi]\to \real$ such that
  \[\dot c (s) = (\cos(\alpha(s)),\sin(\alpha(s))) \]
  then the angle function is embedded in a suitable subspace $N$
  of $L^2(0,2\pi)$ or $W^{1,2}(0,2\pi)$.
  Since the goal is  to obtain a Shape Space 
  representation for Shape Analysis purposes, 
  a  quotient  is then introduced on $N$.

\item Another representation of planar curves for Shape Analysis is found 
  in Younes \cite{Younes:Comp}. In this case
  the angle function is considered mod$(\pi)$.
  This representation is both simple and very powerful at the same time.
  Indeed, it is possible to prove that geodesics do exist
  and to explicitly show examples of geodesics.

\item Metrics of ``geometric'' curves
  (that is, curves up to the choice of parametrization) 
  have been studied by Michor--Mumford
  \cite{Michor-Mumford,MM:hamil_rieman,Michor-Mumford04} and
  Yezzi--Mennucci \cite{YM05:confor_metric,YM:eusipco,YM:metrics04};
  more recently,
  Yezzi--Mennucci--Sundaramoorthi
  \cite{ganesh:new_sobol_activ_contour07,
    ganesh:SAC06,GYM-TPAMI07,
    SYM06:track_with_sobol_activ_contour,
    MYS:PropSAC,
    SYM05:sobol_activ_contour,
    ganesh:SAC}  
  have studied Sobolev--like metrics of curves and shown many good
  properties for applications to Shape Optimization; similar results have also
  been shown independently by Charpiat \emph{et al} \cite{Faugeras:App}.

  \begin{Remark}\label{rem:quot_diffeo}
    In this case, shapes are modeled as immersed parametric curves
    $c:S^1\to\real^N$, for the sake of mathematical analysis;
    a quotient w.r.t the group of possible
    reparametrizations of the curve $c$ (that coincides with the group
    of diffeomorphisms $\mbox{Diff}(S^1)$) is applied afterward to all
    the mathematical structures that are defined (such as the
    manifold of curves, the Riemannian metric, the induced distance,
    etc.).
  \end{Remark}

\end{itemize}

\subsection{Goals}
\label{sec:goal}
We remarked that, in Shape Analysis, shapes are usually considered
``up to rotation, translation and scaling'', but even in Shape
Optimization, to a certain degree, our theory should be independent of
rotation and translation: that is, whatever we do with shapes should
not depend on ``where in the plane'' we do it.

In the rest of the paper we 
will denote by $\calI$  the family  of the nonempty compact sets in
$\real^N$, and we will build many examples of metrics $d$  on $\calI$. We
will always require these metrics to be
\textbf{euclidean invariant}. If $A$ is an euclidean transformation 
of the space (a rigid transformation), then 
\[d(A \Omega_1,A \Omega_2)=d(\Omega_1,\Omega_2)~.\]

What other properties may be interesting for applications?

As mentioned before, a goal of Shape Optimization is to define
\emph{shape metrics},
\emph{shape averages}, \emph{shape principal component analysis},
\emph{shape probabilities}\ldots

For example, if we represent shapes $A_j,j=1\ldots n$
by their \emph{signed distance function} $b_{A_j}$,
then we may define 
\emph{Signed Distance Level Set Averaging}
\begin{equation}
\bar{A}=\Big\{x\,|\, f(x)=0\Big\},
\mbox{ where } f(x)=\frac{1}{N}\sum_{n=1}^N b_{A_n}(x)\label{eq:sflsa}
\end{equation}
A benefit of this definition is that it is easily computable;
a defect is that, if the shapes are far way, then $\bar{A}$
will be empty. Another defect is that this definition is quite
\emph{ad hoc}: it is not coupled with any other structure that
we may wish to add to the Shape Space, such as a metric $d$.
We may then look at the problem in the other direction.

Considering a generic metric space $(M,d)$,  define
the \emph{Distance Based Averaging} 
\footnote{also known as \emph{Karcher mean},
%H. Karcher. Riemannian center of mass and mollifier
%smoothing. Communications on Pure and Applied
%Mathematics, 30:509â541, 1977.
 but it is also sometimes attributed to Fr\'echet, in 1948}
of any given collection $a_1\ldots a_n\in M$, as
a minimum point $\bar{a}$ of the sum of its squared distances:
\begin{equation}
 \bar{a}=\mathop{\arg\min}_a \sum_{j=1}^n d(a,a_j)^2 \label{eq:dba}
\end{equation}

Supposing now that the Shape Space $\calI$ is given a metric $d$,
we can use the abstract definition above to define \emph{shape
averages}; this definition has many advantages. Namely
\begin{itemize}
\item it comes from a minimality criterion, so it is \emph{``optimal''}
  in a certain sense (contrary to the definition \eqref{eq:sflsa}).
\item   If the distance is invariant w.r.t. 
  a group action, then the \emph{shape average} is as well
  (see sec.~\ref{sec:dist_quot}). For example, in the case
  of \emph{geometric curves}, where the distance is independent
  of parametrization, then the  \emph{shape average} will
  be independent of the parametrization of $a_1\ldots a_n$.
\item It coincides with the arithmetic mean in Euclidean spaces; more
  in general, when $\calI$ is a smooth submanifold of a Banach space and
  $a_1\ldots a_n$ are near enough, then  $\bar{a}$ is an approximation of the
  arithmetic mean.
\end{itemize}

In particular, the average of two shapes $A_1,A_2$
is the \emph{midpoint}, that is a shape $A$ such that
\[d(A_1,A)=d(A_2,A) = \frac 1 2 d(A_1,A_2) \]
We are then, however, bound by this result (whose base idea  goes back
to Menger -- see sec. \S4.i.1 % \ref{ASYM-sec:mid-point} 
in \cite{ACM:Asym} for more details)
\begin{Theorem}
  Supposing that $d$ is   complete and \emph{intrinsic} (see sec.~\ref{sec:metric_spaces}), then the
  following facts are equivalent:
  \begin{itemize}
  \item for any two shapes $A_1,A_2$ there 
    is a \emph{midpoint};
  \item for any two shapes $A_1,A_2$ there is a minimal geodesic connecting
    them.
  \end{itemize}
\end{Theorem}
For this reason, we end up studying whether the Shape Space 
admits minimal geodesics
(in theorem \ref{thm:exists!}).

\ammark{COsa altro? aggiungere qualche discussione sullo shooting di geodetiche}

\subsubsection{Tradeoff}
\label{sec:antinomy}
% When dealing with problems in
% that is, finding the 
% there is a known antinomy to address
Unfortunately a tradeoff 
(that is well known in the  Calculus of Variations) arises;
\begin{itemize}
\item on one hand, a Shape Space that is useful for Shape Optimization should
  possibly be equipped with a topology that makes functionals ``regular'',
  so that suitable minimization methods can be used; to this end, the topology
  should have many open sets.

\item On the other hand, to prove existence of average points and 
  of geodesics (that are useful in Shape analysis),
  it is sufficient that certain bounded sets be compact
  (cf. \ref{prp:D_compatto}, \ref{prp:D_completo},  \ref{prp:exist_dgba}):
  to this end, the topology  should have few open sets.
\end{itemize}

We can exemplify this as follows.
\begin{itemize}
\item 
  As we mentioned before, it was shown
  in \cite{ganesh:SAC06,GYM-TPAMI07} that flows for Active Contour
  methods that use a  Sobolev metric are more robust to noise
  and converge faster  than standard  flows.
  To explain the rationale,
  suppose $H,L$ are two metrics, $H$ being stronger than  $L$.
  % so we would expect that the flow induced by $H$ would be less
  %sensitive to noise than the the flow induced by $L$;
  \begin{figure}[htbp]
    \begin{center}
      \includegraphics[width=2in]{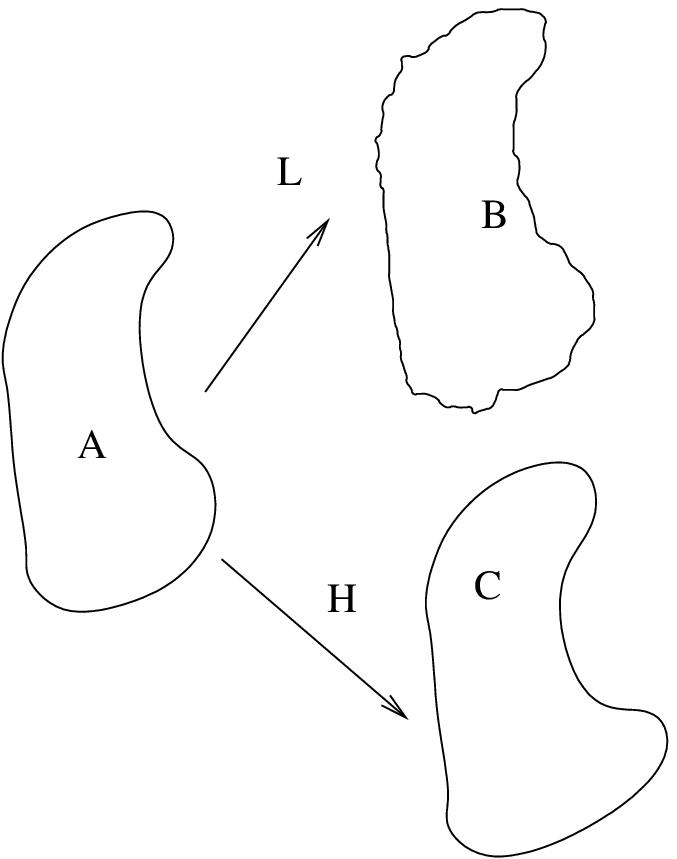}
    \end{center}
    %\caption{example \ref{exa: noncomp simple}}
    %\label{fig: noncomp simple}
  \end{figure}
  When evolving the shape $A$,
  the $L$--flow will move towards a shape $B$ with small scale deformations
  (such as those induced by noise),  since $B$ is nearer to $A$ in the 
  the $L$--induced distance; whereas
  the $H$--flow will move towards
  the shape  $C$ with large scale deformations, 
  since $C$   is nearer to $A$ w.r.t its related distance.
\item Suppose now though that a dataset contains a template shape $A$;
  an algorithm is given a version $B$ of $A$ that was corrupted by noise,
  and different shape $C$, and it has to decide what is the best
  match to $A$. In this case, the weaker metric $L$ would associate
  $A$ to the correct shape $B$, whereas an algorithm
  using the metric $H$ would fail to associate  $A$ to $B$.
\end{itemize}

For all those reasons, it is quite difficult to find a Shape Space
that is suited both for Shape Analysis and for Optimization.

\subsection{Plan of the paper}
The plan of the paper is as follows:
we foremost 
%state and discuss
%our goals, that is to define \emph{geometric} metrics of compact sets; then we
provide base definitions, and we propose some results 
in the theory of metric spaces, in particular when they are
isometrically embedded in Banach spaces. Considering the space $\calI$ of
compact sets
we review the well--known Hausdorff distance, and its properties;
we successively propose a class of metric spaces that are similar to the
Hausdorff distance, while at the same time enjoying some extra
properties that may be useful in applications.

\section{Metric spaces and embeddings in Banach spaces}

\subsection{Metric spaces}
\label{sec:metric_spaces}
We recall some basilar definitions and results
in the abstract theory of metric spaces.
Suppose that $(M,d)$ is a metric space.
We induce from $d$ the length $\len^{d} \gamma$ 
of a continuous curve $\gamma:[\alpha,\beta]\to M$,
by using the total variation
\begin{equation}
  \label{eq: variation as len}
  \len^{d}\gamma \defeq \sup_T \sum_{i=1}^n d\big(\gamma(t_{i-1}),\gamma(t_{i})\big),
\end{equation}
where the sup is carried out over all finite subsets $T=\{t_0,\cdots,t_n\}$ of $[\alpha,\beta]$
and $t_0\le\cdots\le t_n$.

We define the \textbf{induced geodesic distance} $d^g$ by
\begin{equation}\label{eq: geo b}
  d^g(x,y) \defeq \inf_\gamma \len^{d_1} \gamma,
\end{equation}
where the inf is taken in the class of all continuous curves $\gamma$
connecting $x$ to $y$. If the $\inf$ is a minimum, the curve
providing the minimum is called a \emph{geodesic}.
Note that it may be the case that $d^g(x,y)=\infty$ for some choices
of $x,y$. Note also that $d^g\ge d$.

\begin{figure}[htbp]
  \begin{center}
    \includegraphics{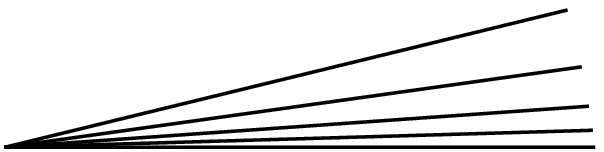}
  \end{center}
  \caption{example \ref{exa: noncomp simple}}
  \label{fig: noncomp simple}
\end{figure}

We may consider the metric space $(M,d^g)$; but
the topology of $(M,d)$ and $(M,d^g)$ may be quite different,
as we see in this example:
\begin{Example}\label{exa: noncomp simple}
\em  Let us consider the  subset $M \defeq \psi(E)$ of $\real^2$,
  where 
  \[E \defeq  [0,1] \times \Big(  \gbrk{ 0 } \cup
  \gbrk{1/ n ~|~ {n>0}}  \Big)\]
  and 
  \[\psi(\rho, \theta) \defeq \rho (\cos(\theta), \sin(\theta)).\]
  (See fig.~\vref{fig: noncomp simple}).
  We associate to $M$ the distance $d$ induced by
  the usual distance of $\real^2$ onto $M$ and $d^g$ for the 
  geodesic distance. Then we have that $(M,d)$ is obviously compact whereas
  $(M,d^g)$ is not: indeed  $x_n\defeq\psi(1,1/n)$ does not admit
  a converging subsequence, since
  $d^g(x_n,y_m)=2$ for all $n\neq m$.
\end{Example}

When $d=d^g$, we will say that the metric space is \textbf{path--metric},
or that $d$ is \textbf{intrinsic}.

Note that the length  $\len^{d^g}$ defined by $d^g$ coincides with
$\len^{d}$, and then $d^g=(d^g)^g$: $d^g$ is always intrinsic.

\ammark{citare riparametrizzazione a arc parameter}

We will use the following proposition:
\begin{Proposition}\label{prp:D_compatto}
  if for a choice of $\rho>0$
  \begin{equation}
    \mathbb D^g(x,\rho)\defeq \{ x ~|~ d^g(x,y)\le \rho\}\label{eq:def:Dg}
  \end{equation}
  is compact  in the $(M,d)$ topology, then
  $x$ and any $y\in \mathbb D^g(x,\rho)$  may be connected by a geodesic.
\end{Proposition}
The proof is simply obtained by the direct method in the Calculus of
Variations (see Thm.~4.24
%\ref{ASYM-Busemann local}
in \cite{ACM:Asym}).

% using Ascoli--Arzel\`a theorem 
% and the lower semi continuity of $\len ^{d}$
% (if  $\lim\gamma_n=\gamma$ pointwise then
% $\liminf _n \len ^{d}\gamma_n\ge \len ^{d}\gamma$).

We also state these simple propositions.
\begin{Proposition}\label{prp:D_completo}
  If for a $x\in M$ for all choices of $\rho>0$
  \[\mathbb D(x,\rho)\defeq \{ x ~|~ d(x,y)\le \rho\}\] is compact
  then $(M,d)$ is complete.
\end{Proposition}

\begin{Proposition}\label{prp:exist_dgba}
  Suppose that $a_1\ldots a_n\in M$ are given;
  a  sufficient condition for the existence of 
  the \emph{Geodesic Distance Based Averaging} $\bar{a}$ of   $a_1\ldots a_n$
  \begin{equation}
    \bar{a}=\argmin_a    \tau(a) \text{ , where }
    \tau(a) \defeq \sum_{j=1}^n d^g(a,a_j)^2 \label{eq:gdba}
  \end{equation}
  is that, defining
  \[ \rho^*=\min_{i=1,\ldots n} \tau(a_i) 
  %~~,~~ i^*=\argmin_{i=1,\ldots n} \tau(a_i) , 
  \]
  we have that $\rho^*<\infty$ and
  that  $\mathbb D^g(a_{1},2\sqrt{\rho^*}+\varepsilon)$ 
  %(see \eqref{eq:def:Dg})
  is compact  in the $(M,d)$ topology,  for $\varepsilon>0$ small.
  \begin{proof}
    Note first that the infimum of $\tau(a)$ 
    is finite, since it does not exceed  $\rho^*$.    Recall that 
    \[d^g(a,a_j) = \inf_{\gamma_j} l_j \] where $l_j$ is the length of a
    Lipschitz curve $\gamma_i$ connecting $a,a_i$. 
    So we can rewrite the problem    \eqref {eq:gdba} as
    \[\inf_{\gamma_1\ldots \gamma_n} \theta (\gamma_1\ldots \gamma_n),\quad
    \text{where } \theta (\gamma_1\ldots \gamma_n)\defeq    \sum_{j=1}^n (l_j)^2\]
    where the infimum is computed on all choices of Lipschitz curves
    $\gamma_1\ldots \gamma_n$ of length $l_1\ldots l_n$ connecting
    $a_i$ to a common point $x\in M$; for simplicity we represent them 
    as $\gamma_i:[0,l_i]\to M$ parametrized by arc parameter.
    By    triangular inequality
    \[d^g(a_i,\gamma_{j}(t)) \le d^g(a_i,x) + d^g(x,\gamma_{j}(t))
    \le l_i+ l_j \]
    Let then $\gamma_{i,k}$ be a sequence of choices
    that converges to the infimum:
    \[    \theta (\gamma_{1,k}\ldots \gamma_{n,k})\to_k 
    \inf_{\gamma_1\ldots \gamma_n} \theta (\gamma_1\ldots \gamma_n)    \]
    so for large $k$, 
    \[\theta (\gamma_{1,k}\ldots \gamma_{n,k})\le \rho^*+\varepsilon\]
    but then in particular $l_{i,k}\le \sqrt{\rho^*+\varepsilon}$
    hence
    \[d^g(a_{1},\gamma_{j}(t)) \le 2  \sqrt{\rho^*+\varepsilon} \]
    so all the curves are contained in a compact set.
    By Ascoli--Arzel\`a theorem,
    we can then extract a uniformly convergent subsequence, and use the 
    fact that the length is lower semi continuous.    
  \end{proof}
\end{Proposition}

A similarly proposition can be stated for $d$:
\begin{Proposition}\label{prp:exist_dba}
  Suppose that $a_1\ldots a_n\in M$ are given;  let
  \begin{equation}\label{eqn:rho_avg}
    \rho^*=\min_i \sum_{j=1}^n d(a_i,a_j)^2
  \end{equation}
  and  $i^*$ the index that achieves the above minimum:
  suppose that  $\mathbb D(a_{i^*},\sqrt{\rho^*}+\varepsilon)$ is compact
  for $\varepsilon>0$ small:  then there exists a point  $\bar{a}$
  that is  the \emph{ Distance Based Averaging} of 
  $a_1\ldots a_n$, as defined in \eqref{eq:dba}.
  \ammark{penso che si possa fare anche con $\varepsilon=0$}
\end{Proposition}

\subsubsection{Distances, quotients and groups}
\label{sec:dist_quot}
Let $d_M(x,y)$ be a distance on 
a space $M$, and $G$ a group acting on $M$;
a distance $d_B$ may be defined on $B=M/G$ by 
\[d_B ([x],[y])= \inf_{x\in [x],y\in [y]} d_M(x, y)
= \inf_{g,h\in G} d_M(g x, h y)\]
that is the lowest distance between two orbits; 
% given $\bar x,\bar y\in B$ equivalence classes of shapes,
we write $d_B (x,y)$ for simplicity.

If $d_M$ is \textbf{ invariant w.r.t. $G$}, i.e.
\[d_M(g x, g y) = d_M(x, y) \quad \forall g\in G\]
then
\begin{equation}
d_B (x,y) = \inf_{g\in G} d_M(g x,  y)\label{eq:eq:inv_quot}
\end{equation}
It is easy to see that $d_B$ satisfies the triangular inequality;
but it may be the case that $d_B(x,y) =0$ even when $x\neq y$.
%that is, $d_B$ is a semidistance. 
We state a simple sufficient condition
\begin{Lemma}\label{lemma:quot}
   If the orbits are compact, then $d_B$ is a distance.
\end{Lemma}

When studying metrics $d$ on a Shape Space $\calI$,
the quotient is particularly useful in at least  two cases:
\begin{itemize}
\item when we want to pass from a \emph{preshape space}
  \footnote{cf. remark \ref{rem:preshape}} to a \emph{shape space}: in
  this case, $G$ is the Euclidean group of
  rotations and translation (and sometimes of  scaling);
\item when the representation is redundant:
  for example, in remark \ref{rem:quot_diffeo} we would set $G=\mbox{Diff}(S^1)$.
\end{itemize}

\ammark{Can we say that $(d_B)^g$ is equal to  the projection to the quotient
of $(d_M)^g$?}

%We will use the above both to define shape spaces and to pass
%from preshapes to shapes.

\subsection{Embeddings in Banach  spaces}
\label{sec:emb}
In most of what follows, we will be able to identify $M$ (using an
isometry $i$) with a subset $N$ of a Banach space $E$.
We remark that an isometry is a map $i$ such that
$d(x,y)=\|i(x)-i(y)\|$ (and this should not be confused with
the concept of  isometrical embedding of Riemannian manifolds).

\subsubsection{Radon-Nikodym property}

The following result
from \cite{Ambrosio-Kircheim00} will come handy:
\begin{Theorem}\label{thm:len_d_B}
  Suppose that $E$  is the dual of a separable Banach space.
  Let $\gamma:[a,b]\to E$ be a Lipschitz curve; then, by thm.~8.1 in 
  \cite{Ambrosio-Kircheim00},  for almost
  all $t$ there exists the derivative $\dot \gamma(t)\in B$
  that is defined as
  \begin{equation}
    \dot \gamma(t) \defeq 
    \text{w-}\lim_{\tau\to 0} \frac{\gamma(t+\tau)-\gamma(t)}\tau\label{eq:w_der}
  \end{equation}
  where the limit is done according to the weak-* topology;
  and moreover, 
  \begin{equation}
    \|\dot \gamma(t)\|=
    \lim_{\tau\to 0} \left\|\frac{\gamma(t+\tau)-\gamma(t)}\tau\right\|
    \label{eq:w_der_norm}
  \end{equation}
  so $\|\dot \gamma(t)\|$ coincides with the 
  \emph{metric derivative}, that is studied in 
  \cite{Ambrosio-Tilli}.

  There follows easily (by applying scalar products to \eqref{eq:w_der}) that
  \begin{equation}
   \gamma(b)-\gamma(a) = \int_a^b \dot \gamma(t) dt\label{eq:gamma_a_b}
 \end{equation}
  and
  \begin{equation}
    \len^d\gamma = \int_a^b \|\dot \gamma(t)\| dt
    \label{eq:len_b_dot_gamma}
  \end{equation}
  (this last by thm.~4.1.1 in \cite{Ambrosio-Tilli}).
\end{Theorem}
It is common to say that $E$ enjoys the \emph{Radon-Nikodym Property}, when
the limit in \eqref{eq:w_der} exists in the strong sense, and for almost all
$t$.
Note that the \emph{Radon-Nikodym Property} does not hold in general:
consider the map $t\mapsto \uno_{[t,t+1]}$ in $L^1(\real)$, whose
derivative should be $t\mapsto \delta_{t+1}-\delta_t$.

\smallskip

We now recall this basilar definition:
\begin{Definition}
  a Banach space $E$ is  \textbf{uniformly convex} if
  $\forall \e>0 \exists \delta>0$,
  \[ \forall x,y\in E, \|x\|\le 1,\|y\|\le 1, \|x-y\| \ge \e
  \Longrightarrow \|(x+y)/2\|< (1-\delta)~.\] 
\end{Definition}
Examples of uniformly
convex Banach spaces include $L^p(\Omega,\mathcal A,\mu)$
for $p\in (1,\infty)$.
\quad
Uniformly  convex Banach spaces have many interesting 
properties: for example, they are reflexive
(Milman Theorem,  III.29 in \cite{Brezis}); moreover, if
$x_n\to x$ in weak sense and $\limsup \|x_n\|\le \|x\|$ then
$x_n\to x$ in the strong sense (prp III.30 in \cite{Brezis}).

So we obtain  a sufficient condition:
\begin{Corollary}\label{Radon-Nikodym}
  if $E$ is uniformly convex and separable, then it enjoys the
  \emph{Radon-Nikodym Property} (indeed
  % by the cited prp III.30 in \cite{Brezis},
  eqn.~\eqref{eq:w_der} and eqn.~\eqref{eq:w_der_norm} imply that the
  limit in \eqref{eq:w_der} is valid also in the strong sense).
\end{Corollary}

\subsubsection{Embeddings in uniformly convex Banach spaces}
If  $E$ is  uniformly convex then in particular
the closed ball $\{x~|~\|x\|\le 1\}$ is strictly convex; this
has a curious implication.
\begin{Lemma}
  Suppose the closed balls in $E$ are strictly convex.
  Consider $E$ as a metric space, with distance $d_E(x,y)=\|x-y\|$.
  The segment connecting $x,y\in E$ is the unique minimal geodesic
  (up to reparametrization).
  \begin{proof}
    We will prove that, for $x,y$, for any
    minimal geodesic $\gamma:[0,1]\to M$ connecting $x$ to $y$,
    if $\gamma$ is reparametrized to arc parameter then
    $\gamma(1/2)=(x+y)/2$; iterating this reasoning
    with finer subdivision  we obtain that $\gamma(t)=(t x+ (1-t)y)$.

    With no loss of generality, up to translation
    and scaling, suppose  $y=-x$  and $\|x\|=1$.
    The segment $t\mapsto tx$ is a geodesic for $t\in [-1,1]$,
    by the theorem \ref{thm:len_d_B}, and its length is 2.
    Suppose now that $\gamma:[-1,1]\to M$ is another geodesic: then
    $\len \gamma=2$, and, up to reparametrization,  $\|\dot\gamma\|=1$
    at almost all points; in particular, setting $z=\gamma(0)$,
    $\|z-y\|\le 1$ and $\|x-z\|\le 1$;
    but then, by triangular inequality, 
    $ \|z+x\| = \|x-z\| = 1$.
    Suppose that $z\neq 0$; then $\|(z+x)-(x-z)\|> 0$;
    by strict convexity, though, this implies that
    $\| ((z+x)+(x-z))/2 \|= \|x\|<1$, and this is a contradiction.
  \end{proof}
\end{Lemma}

\begin{Theorem}\label{thm:geodetiche_tese}
  Suppose that 
  $(M,d)$ is a complete space, and that
  $i:M\to E$ is an isometrical immersion in a
  uniformly convex Banach space $E$. If, given $x,y\in M$,
  $d(x,y)=d^g(x,y)$, then the segment connecting $i(x),i(y)$
  is all contained in $i(M)$.

  In particular, if $(M,d)$ is  path-metric then 
  $i(M)$ is convex, and then any two points in $M$ can
  be joined by a unique minimal geodesic   (unique up to reparametrization).
  \begin{proof}
    Note that  $i(M)$ is complete, and then it is closed in $E$.
    We will prove that, for any $x,y\in i(M)$, $(x+y)/2\in i(M)$;
    we can then iterate this idea to further subdivide, and 
    since $i(M)$ is closed then this proves the whole
    segment connecting $x,y$ is in $i(M)$; for the above lemma,
    the segment is the unique minimal geodesic.

    We now fix $x,y\in i(M)$:  there must be  paths
    $\gamma_n:[-1,1]\to i(M)$ connecting $x$ to $y$ with 
    length $\len(\gamma)< L_n\defeq \|x-y\|+2/n$.

    As in the lemma before, we suppose for simplicity that
    $y=-x$  and $\|x\|=1$ (so $L_n=2+2/n$); and we
    reparametrize so that   $\|\dot \gamma_n\|=1+1/n$: hence
    setting $z_n=\gamma_n(0)$
    \[\|z_n+x\|\le 1+1/n ~~,~~\|x-z_n\|\le 1+1/n~.\]
    and then by triangle inequality $\|z_n+x\|\to 1$,
    $\|z_n-x\|\to 1$.
    Setting \[w_n=(z_n+x) /\|z_n+x\|~~,~~v_n=(x-z_n) /\|z_n-x\|\] 
    we can prove that $ \| (w_n+v_n)/2\| \to 1$ hence
    by the uniform convexity of $E$ we obtain that  $w_n-v_n\to 0$ 
    and then $z_n\to 0$. Since $z_n\in i(M)$
    and $i(M)$ is closed then $0\in i(M)$.
  \end{proof}
\end{Theorem}
% So, when a complete metric space is isometrically embedded
% in a uniformly  convex Banach space $E$, either $i(M)$ is convex,
% or $(M,d)$ is not path metric.

The above is a ``rigidity theorem'', in that it restricts the class of
metric spaces that can be isometrically embedded in a uniformly convex
Banach space~$E$.
\begin{Corollary}
  a complete compact finite dimensional Riemannian manifold $M$ cannot
  be isometrically embedded in a uniformly convex Banach space $E$:
  indeed in this space $M$ there are two points
  that can be joined by more than one minimal geodesic.
  % , and  this is impossible since $i(M)$ would be convex.
\end{Corollary}

When $E$ is not uniformly convex, on the other hand, strange behaviours
arise.
\begin{Remark}\label{rem:brutto}
  Let $L^\infty=L^\infty(\Omega,\mathcal A,\mu)$ and suppose $\Omega$
  is not an atom of $\mu$, that is, suppose the dimension of 
  $L^\infty$ is greater than 1. Given generically $f,g\in L^\infty$, there is an
  uncountable number of minimal geodesics connecting them.
  \begin{proof}
    We can assume without loss of generality that $g=0$ and that
    $\|f\|=1$. Let $A=\{|f|=1\}$.  \quad We will prove that if there
    is only one geodesic then $|f|=\uno_{A}$.  \quad Indeed if
    $|f|\neq \uno_{A}$ then $\mu\{|f|<1\}>0$. Let $0<t<1$ be such that
    $\mu\{|f|<t\}>0$; obviously $\mu\{|f|\ge t\}>0$ since $\|f\|=1$;
    let $A'=\{|f|\ge t\}$ and $A''=\{|f|< t\}$.  \quad Given any
    diffeomorphism $b:[0,1]\to[0,1]$ with  $b'(s)\le 1/t$,
    \[ \gamma(t)\defeq t f \uno_{A'} + b(t) f \uno_{A''} \]
    is a geodesic. Indeed its derivative is
    \[ \gamma'(t)\defeq f \uno_{A'} + b'(t) f \uno_{A''} \]
    and $\|\gamma'(t)\|=1$ by construction.

    The family of $f$ s.t. $|f|=\lambda\uno_{A}$ is closed and
    has empty interior.
  \end{proof}
\end{Remark}

The idea of \emph{isometrical embedding} is quite powerful:
indeed any separable metric space
may be isometrically embedded in $\ell^\infty$ (that is the dual of the
separable space $\ell^1$): so the breadth of application of the 
 theorem \ref{thm:len_d_B} is general, and is at the basis
of many results in \cite{Ambrosio-Kircheim00}.
But  the embedding in $\ell^\infty$ that is
studied in  \cite{Ambrosio-Kircheim00} is  not suited
for our practical applications:
\begin{itemize}
\item it would not respect the geometric properties of the space
  (as we discussed in sec.~\ref{sec:goal})
\item  it would be too difficult
  to find a satisfactory notion of \emph{``shooting of minimal
    geodesics''} using this embedding.
\end{itemize}

For all above reasons, we will consider \emph{isometrical
  embeddings} in this paper as well but we will
(for the most interesting applications)
use an explicitly chosen embedding in uniformly convex Banach spaces.

\subsection{Definitions}
We introduce some definitions that will be used in the rest of the paper

We will write $s\vee t = \max\{s,t\}$ and 
$s^+=\max\{s,0\}$, when $s,t\in\real$.

We will write $B(x,r)$ or $B_r(x)$ for the open ball of center 
$x$ and radius $r>0$ in $\real^N$; we will shortly write $B_r$ for $B_r(0)$.
\quad
Similarly $D_r(x)$ will be the closed ball of center 
$x$ and radius $r>0$ in $\real^N$, and $D_r=D_r(0)$.

We define the \textbf{fattened set} to be
\[ %\fat{A}{r} \defeq 
A + D_r = \{ x+y ~|~ x\in A , |y| \le r \} =
\bigcup_{x\in A} D_r(x)= \{ y ~|~ u_A(y) \le r \}.\]
This fattened set is always closed,
(since the distance function  $u_A(x)$, that was defined in  \eqref{eq:distfun}, is continuous).

We will say that a family $A_{i\in I}$ of sets in $\real^N$ is 
\emph{equibounded} if there is a $R>0$ such that $A_i\subset D_R$ for all $i$.

We denote by $\mathcal{L}^N$ the $N$ dimensional Lebesgue measure,
and $\omega_N\defeq \mathcal{L}^N(B_1)$;
we write shortly $\int_A f(x)dx$ for the Lebesgue integral.

\subsection{Hausdorff distance}

A fundamental example of metric on $\calI$ 
is the \textbf{Hausdorff distance}
\[ d_H(\Omega,\Omega') \defeq \inf \{ \delta > 0 ~|~ 
\Omega'\subset (\Omega + D_\delta) ,\Omega\subset (\Omega' + D_\delta)\} ~~. \]
It is not difficult to verify that
\begin{equation}
  d_H(\Omega,\Omega') = \sup _{x\in\real^N} |u_\Omega(x) - u_{\Omega'}(x) | ~~,
  \label{eq:dh_da_sup}
\end{equation}
see for example  Thm.~2.2 in ch.~4 in Delfour--Zolesio \cite{Delfour2001}.

This metric enjoys many important properties.
\begin{Theorem}\label{thm:dH}
  The metric space $(\calI,d_H)$ satisfies:
  \begin{itemize}
  \item given $r>0$, the family of $r$-bounded compact sets
    \[\{ \Omega\in\calI ~|~ \Omega\subset D_r \}\]
    is compact; in particular, the set 
    \[\mathbb D\defeq \{ \Omega ~|~ d_H(\Omega,\Omega')\le \rho\}\] is compact;
  \item $(\calI,d_H)$ is path--metric
    (that is, $d_H=(d_H)^g$)
  \item consequently, by Prp.~\ref{prp:D_compatto}
    any two $\Omega,\Omega'\in \calI$ may be joined by a minimal geodesics;
  \item and moreover, by Prp.~\ref{prp:D_completo},
    $(\calI,d_H)$ is complete.
  \end{itemize}
\end{Theorem}

The first statement is a well known property of
the Hausdorff distance, see e.g.  \cite{Delfour2001} pag. 194.
By exploiting the characterization \eqref{eq:dh_da_sup},
it also follows from a diagonal/compactness argument
and the following rigidity property:
\begin{Lemma}\label{lem:conv_punt}
  Let $\Omega_n$ be closed sets, and
  suppose that $\lim_n u_{\Omega_n}(x) = f(x)$
  for all $x$ in a dense subset $D$ of $\real^N$.
  Then there is a closed set $\Omega$ such that 
  $u_\Omega(x) = f(x)$ for all $x\in D$, and
  $u_{\Omega_n}\to u_\Omega$ uniformly on compact sets;
  moreover  if (and only if) $\Omega_n$ is equibounded then
  $u_{\Omega_n}\to u_\Omega$ uniformly.
  \begin{proof} 
    The proof may follow from the theory of Viscosity Solutions:
    it is well known, indeed, that $u_\Omega$ is the unique solution
    to a properly defined \emph{Eikonal equation}; and that 
    viscosity solutions do enjoy the required rigidity property.

    We propose here instead a direct proof.
    We set $u_n \defeq u_{\Omega_n}$;
    it is easily proved that $u_n$ is 1-Lipschitz, that is
    \begin{equation}
      \label{eq:1lip}
      |u_n(x)-u_n(y)|\le |x-y| \quad \forall x,y
    \end{equation}
    so  passing to the limit in the above \eqref{eq:1lip},
    we obtain
    \begin{equation}
      \label{eq:1lip_f}
      |f(x)-f(y)|\le |x-y| \quad \forall x,y\in D
    \end{equation}
    and then there is an unique  extension
    of $f$ to a positive   function $g:\real^N\to\real$
    that is again 1-Lipschitz, that is, 
    \begin{equation}
      \label{eq:1lip_g}
      |g(x)-g(y)|\le |x-y| \quad \forall x,y~~.
    \end{equation}

    It is easy to prove that $u_n(x)\to g(x)$ for all $x$,
    and actually (by imitating the proof of Ascoli--Arzel\`a theorem) that
    $u_n\to g$ uniformly on compact sets.
    
    Let $\Omega=\{g=0\}$; to conclude the proof, we need to prove
    that $g=u_\Omega$. To this end, we first prove that
    $g\ge u_\Omega$: indeed, fixing $x$,
    $u_n(x)=|x-y_n|$ for at least one point $y_n\in\Omega_n$;
    since  $u_n(x)\to g(x)$, then the sequence $\{y_n\}$ is bounded,
    so (up to a subsequence $n_k$) it converges to a point $y$; since the family
    $u_n$ is 1-Lipschitz and $u_n(y_n)=0$ then
    $g(y)=0$, that is $y\in \Omega$: hence 
    \[g(x)=\lim_k u_{n_k}(x)= \lim_k|y_{n_k}-x| = |y-x|\ge u_\Omega(x)~.\]
    Conversely, let $y\in \Omega$ be such that $u_\Omega(x)=|x-y|$; then
    by \eqref{eq:1lip_g}
    $g(x)\le g(y) + |x-y|=|x-y|=u_\Omega(x)$.

    To conclude, supposing that   $\Omega_n$ is equibounded,
    then choosing $R>0$ such that $\Omega_n\subset D_R$, we know that 
    $u_{\Omega_n}\to u_\Omega$ uniformly on $D_R$, so given $\e>0$ 
    for $n$ large $|u_n-u|<\e$  in $D_R$ and then
    \[u_n(x)=\inf_{y\in D_R}(|x-y|+u_n(y)) \le 
    \inf_{y\in D_R}(|x-y|+u(y)+\e) = u(x)+\e\]
    where the first and last equalities are due to the Dynamical Programming
    principle; and similarly we obtain that $u(x)\le u_n(x)+\e$.
    \quad
    The ``only if'' part follows from \eqref{eq:dh_da_sup}.
  \end{proof}
\end{Lemma}

To prove the above second property
in \ref{thm:dH}, we may use the first property
and the following \emph{Menger convexity} result
\begin{Proposition} Let $A,B \in \calI$ be two compact sets, then for all $\lambda \in [0,1]$there exists a compact set $C$ such that
  \begin{itemize}
  \item $d_H(A,C) = \lambda d_H(A,B)$,
  \item $d_H(B,C) = (1 - \lambda) d_H(A,B)$.
  \end{itemize}
  \begin{proof}
    We write $[A]_r=A+D_r$ for the fattened set.
    Let $\mu = d_H(A,B)$. We consider the set 
    $$
    C \defeq \left\{ z | \exists x \in A,y \in B, \, |x-y| \leq \mu,
      |x-z| \leq \lambda \mu, \, |y-z| \leq (1-\lambda) \mu \right\}.
    $$
    We prove that the set $C$ has the properties we need. In particular it is enough to prove only the first one because of the symmetry in the two conditions.
    If $x \in A$ then there exists $y \in B$ such that $|x-y| \leq \mu$, then
% the geodesic condition on the metric $d$ ensure the existence of 
 the a point $z=(1-\lambda )x+\lambda y$ satisfies
    $$|x-z| \leq \lambda \mu, \, |y-z| \leq (1-\lambda) \mu.$$
    Such a $z$ must be an element of $C$ and so we found an element of $C$ with distance less or equal to $\lambda \mu$. This means that $x \in [C]_{\lambda \mu}$ and it is true for all $x\in A$ so 
    $A \subset [C]_{\lambda \mu}$.

    Let's take now $z \in C$. From the definition of the elements of $C$ we have that there must exists $x \in A$ and $y \in B$ such that $|x-z| \leq \lambda \mu, \, |y-z| \leq (1-\lambda) \mu$. This means that $z \in [A]_{\lambda \mu}$. This is true for all $z \in C$ so $C \in [A]_{\lambda \mu}$.

    To finish the proof we have to show that the set $C$ is compact. It is clearly bounded because it is contained by $[A]_{\lambda \mu}$. We have to show that it is closed. Suppose we have a sequence $\{z_k \}_k \subset C$ such that $z_k \to z$. Then for each $z_k$ we can find two elements $x_k \in A, y_k \in B$ with the properties:
    $$|x_k-z_k| \leq \lambda \mu, \, |y_k-z_k| \leq (1-\lambda) \mu.$$
    The sets $A$ and $B$ are compacts so we can chose a subsequence (for simplicity we use the same index $k$) such that $x_k \to x \in A$ and $y_k \to y \in B$. It is obvious to see that the points $x,y,z$ satisfy the following inequalities:
    $$|x-y| \leq \mu, \ |x-z| \leq \lambda \mu, \ |y-z| \leq (1-\lambda) \mu.$$
    This means that $z$ is an element of $C$ and this concludes the proof.
  \end{proof}
\end{Proposition}

Unfortunately  $(\calI,d_H)$ is quite ``unsmooth'', as shown by this example
(that is similar to \ref{rem:brutto} -- and for a reason!).
\begin{Example}\label{exa:geodetica_dH}
  There are choices of $\Omega,\Omega'\in \calI$ that may be joined by
  an uncountable number of geodesics.

  In fact we can consider this simple example:

  \begin{minipage}{0.4\linewidth}
    \begin{center}
      $A = \left\{ %(x,y) \in M |
        x = 0, 0 \leq y \leq 2 \right\}$
      \\ \medskip
    $B = \left\{ %(x,y) \in M | 
      x = 2, 0 \leq y \leq 1 \right\}$
      \\ \medskip
      $C_t = \left\{ %(x,y) \in M |
        x = 1, 0 \le y \leq \frac 3 2 \right\} \cup
      \left\{ %(x,y) \in M |
        y=0, 1 \leq x \leq t \right\} $
      \\ \medskip
      with $1\le t \le \sqrt 5/2 $;
    \end{center}
  \end{minipage}
  \begin{minipage}{0.55\linewidth}
    \vskip 1pt      \includegraphics{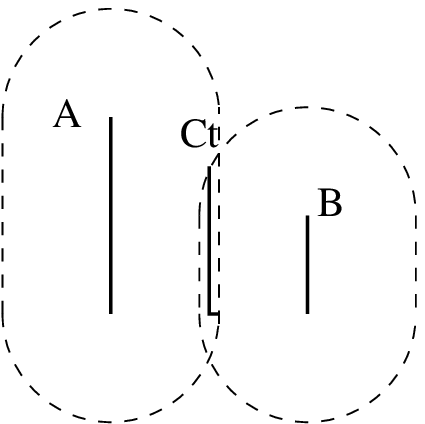}
  \end{minipage}
  and in the picture we represent (dashed) the fattened sets
  $A + B_{\sqrt 5/2}$ and $B + B_{\sqrt 5/2}$.
  Note that  $d_H(A,B) = \sqrt 5$ while
  $d_H(A,C_t) = d_H(B,C_t) = \sqrt 5 /2$: so $C_t$ are all midpoints
  that are on different geodesics between $A$ and $B$.
\end{Example}

We conclude with a family of nice properties.
\begin{Proposition}\rm
  \label{prp:properties_in_dH}
  \begin{enumerate}
  \item    The fattening map $\lambda\mapsto A+D_\lambda$ is  Lipschitz
    (of constant one).
    %since  $d_H(A+D_r, A+D_\lambda)=|r-\lambda$
  \item Given $\lambda>0$, the
    \emph{``fattened area map''}  $L_\lambda(A)\defeq\mathcal{L}^N(A +D_\lambda)$ is  continuous.
  \item Consequently, the area map $L(A)\defeq\mathcal{L}^N(A)$ is upper
    semi continuous.
  \item    Let $\#:\calI\to \N\cup\infty$ be the number $\#\Omega$
    of connected components of a closed set $\Omega$.
    Then $\#$ is lower semi continuous in the metric space $(\calI,d_H)$. 
    
    As a corollary, the family of connected compact sets is a closed
    family in  $(\calI,d_H)$.
    \ammark{\\QUESTE NON LE HA PIU' MOSTRATE NESSUNO
  \item let $\mathbb{G}=\real^N\times \mathcal{S}(N)$ be the group
    that is the (semi)direct product of translations and invertible
    transformations of $\real^N$; then the action on $\calI$ is
    locally Lipschitz, that is, given $A$ compact and $f\in \mathbb{G}$,
    for small $\varepsilon>0$ there exists $\lambda>0$ such that
    \[|f-g|<\varepsilon \Longrightarrow  d_H(f(A),g(A))\le \lambda |f-g| \]
   \item\label{item:removing_in_dH}
     \emph{(removing)}
     Let $x\in A$, then the removing motion
     $\lambda\mapsto A\setminus B_\lambda$ is Lipschitz.
     (See also \ref{exa:removing_in_Lp})
}
  \end{enumerate}
  \begin{proof}
    \begin{enumerate}
    \item Obvious.
    \item if $A_n\to A$ then for fixed $\varepsilon>0$ and
      definitively in $n$,
      \[ A_n\subset A+D_\varepsilon \quad,\quad
      A\subset A_n+D_\varepsilon     \]
      and then
      \[ A_n+D_\lambda\subset A+D_{\varepsilon+\lambda} \quad,\quad
      A+D_{\lambda-\varepsilon}\subset A_n+D_\lambda\]
      passing to Lebesgue measures,
      \[    {\mathcal L}^N(A+D_{\lambda-\varepsilon})\le
      \liminf_n{\mathcal L}^N(A_n+D_\lambda)\le
      \limsup_n{\mathcal L}^N(A_n+D_\lambda)\le
      {\mathcal L}^N( A+D_{\varepsilon+\lambda}) \]
      and we let $\varepsilon\to 0$.
    \item  Since it is the pointwise limit  $L_\lambda(A)\downarrow L(A)$ for $\lambda\to 0$.
    \item     See Thm.~2.3 in ch.~4 in \cite{Delfour2001}. % pag. 158-159.
    \end{enumerate}
  \end{proof}
\end{Proposition}

\section{$L^p$--like metrics of shapes}
The definition of the  Hausdorff distance 
by eqn.~\eqref{eq:dh_da_sup} leads us back to
the paradigm of \emph{representation/embedding};
but in this case it is unfortunately not precise, since the Banach
metric that we use, namely 
\[\|f\|=\|f\|_\infty \defeq \sup_x |f(x)| \]
is usually associated to the spaces $C_b(\real^N)$ of
bounded functions --- whereas the distance function
$u_A$ is not bounded!
What follows is a simple yet effective workaround.

\bigskip
\begin{Definition}\label{H1}\em
  We fix $p\in[1,\infty]$; we fix a function
  $\varphi:[0,\infty)\to(0,\infty)$  monotonically
  decreasing and of class $C^1$, such that
  \begin{equation}
    \varphi(|x|)\in L^p(\real^N). \label{eq:phi_Lp}
  \end{equation}
  Note that, for $p<\infty$, the above is equivalent to asking that
  \begin{equation}
    \int_0^\infty t^{N-1} \varphi(t)^p dt<\infty \label{eq:phi_Lp_bis}
  \end{equation}
  and it implies that $\lim_{t\to\infty} \varphi(t)=0$; for $p=\infty$ we instead
  ask that $\lim_{t\to\infty} \varphi(t)=0$ as an extra hypothesis.
\end{Definition}

An example of such a function is   $\varphi(t)= \exp(-t)$,
or $\varphi=(1+t)^{-(N+1)/p}$.

We will often write $v_A=\varphi\circ u_A$ for simplicity.

\begin{Lemma} \label{lemma:compact_in_Lp}
  Let $\Omega \subset \real^N$ be closed and non empty; 
  suppose $p<\infty$; then 
  \begin{itemize}
  \item[(a)] $v_\Omega \in L^p(\real^N)$ if and only if
  \item[(b)] $\Omega$ is bounded (and then $\Omega$ is compact).
  \end{itemize}
  \begin{proof}
    We first prove that $(a) \implies
    (b)$ by contradiction.  Let us assume that $\Omega$ is
    unbounded. Then there exists a sequence $\set{x_k}\subset \Omega$ such
    that $|x_k|\to \infty$ and $d(x_k,x_q) >2$ for all $k,q \in \N,k
    \not= q$.  The     sequence of sets $B_1(x_k)$
    is disjoint. It is easy to see that $v_\Omega(x) > \varphi(1)$ for
    $x\in\bigcup_k B_1(x_k)$,    and then $v \not\in L^p$.

    Then we prove that $(b) \implies (a)$.  If $\Omega$ is bounded we can
    find a ball $B_R$ such that $\Omega \subset B_R$. Then
    easily we have $u_\Omega \geq u_B \implies v_\Omega \leq v_B$, but $v_B\in
    L^p$ (as is easily proved by $v_B(x)=\varphi( (|x|-R)^+)$ 
    and by \eqref{eq:phi_Lp_bis}) and then also $v_\Omega \in L^p$.
  \end{proof}
\end{Lemma}

\smallskip

\begin{Definition}
  Given $A,B\in\calI$, we define
  \[d_{p,\varphi}(A,B) \defeq \|\varphi(u_A) - \varphi(u_B)
  \|_{L^p(\real^N)} \]
\end{Definition}
By the above lemma, this distance is finite.
We will often write $d$ for $d_{p,\varphi}$
in the following, for simplicity.

\smallskip

The above distance 
is obtained by the  \emph{representation} of a shape $A$
as $v_A$, combined with the 
\emph{embedding}  of $v_A$ in $L^p(\real^N)$.
For this reason, we may identify our shape space with
\begin{equation}
  N_c \defeq \set{v_\Omega~ |~ \Omega \in \calI} 
\end{equation}
that is a subspace of $L^p$.
\begin{Remark}\label{iso}
  By the definition of $d$, the map $\Omega \mapsto v_\Omega$ is an
  isometrical embedding of $\calI$ inside $L^p$, and the image is
  $N_c$; $N_c$ is a closed subset of $L^p$, by the completeness result
  \ref{prp:d_complete} that we will prove in the following.
\end{Remark}
We will exploit this embedding
in the following, as in \S\ref{sec_riemannian}.  

\medskip

It is immediate to verify that $d_{p,\varphi}$
satisfies these properties.
\begin{itemize}
\item The embedding $A\mapsto v_A$ is injective: if $v_A\tilde =v_B$
  then $u_A\tilde =u_B$ (since $\varphi$ is monotonically decreasing,
  and so it is injective); but, by lemma \ref{lem:conv_punt}, this
  implies that $u_A =u_B$ and then $A=B$; consequently, for all
  $A,B\in\calI$, $d_{p,\varphi}(A,B)=0$ iff $A=B$.
\item $d_{p,\varphi}$ is euclidean  invariant,
  as we requested in sec.~\ref{sec:goal}.
\item  
  \begin{equation} \label{eq:d1bound}
    d_{p,\varphi}(\Omega_1,\Omega_2) < \| v_{\Omega_1}\|_{L^p} + \| v_{\Omega_2} \|_{L^p}.
  \end{equation}
  for $p<\infty$, and
  \[d_{\infty,\varphi}(\Omega_1,\Omega_2) < \varphi(0).\]
  \begin{proof}
    When $p<\infty$, by the Minkowski inequality we have that
    $d(\Omega_1,\Omega_2) \le \| v_{\Omega_1}\|_{L^p} + \|
    v_{\Omega_2} \|_{L^p}$; moreover equality would hold only if
    $v_{\Omega_1}=-v_{\Omega_2}$ and this is impossible;
    when $p=\infty$ we use the fact that $\varphi>0$.
  \end{proof}
\item \emph{(Separation at infinity)}
  given two bounded sets $\Omega_1,\Omega_2$ we have
  \begin{equation}
    \lim_{|\tau| \to \infty} d_{p,\varphi}(\Omega_1,\Omega_2 + \tau) = \|
    v_{\Omega_1} \|_{L^p} + \| v_{\Omega_2} \|_{L^p} ~;\label{eq:separ}
  \end{equation}
  for $p<\infty$, and
  \begin{equation}
    \lim_{|\tau| \to \infty} d_{\infty,\varphi}(\Omega_1,\Omega_2 + \tau) = 
    \varphi(0) ~;\label{eq:separ_infty}
  \end{equation}
  \begin{proof}
    For the case $p<\infty$ this comes from a general result for $L^p$
    functions; for $p=\infty$ it derives from the hypothesis
    $\lim_{t\to\infty}\varphi(t)=0$.
  \end{proof}
\item \emph{(Scaling)} If $p<\infty$ and $\lambda>0$ is a rescaling of the space, then
  the rescaled distance may be expressed as
  % \begin{Remark} \label{rmk:nobuono}
  \begin{equation}
    d_{p,\varphi}(\lambda \Omega_1, \lambda \Omega_2)
    =\lambda^{N/p}  d_{p,\tilde \varphi}(\Omega_1,  \Omega_2)\label{eq:scaling}
  \end{equation}
  where  $\tilde \varphi(r)=\varphi(\lambda r)$; indeed
  \begin{align}
    d_{p,\varphi}(\lambda \Omega_1, \lambda \Omega_2)^p
    &=
    \int \lbrk{ v_{\lambda \Omega_1}(x) - v_{\lambda \Omega_2}(x)}^p d\,x \label{eq:zrfsA}\\
    &=
    \lambda^N \int \lbrk{ v_{\lambda \Omega_1}(\lambda z) - v_{\lambda \Omega_2}(\lambda z) }^p d\,z \label{eq:zrfsB}\\ 
    &=
    \lambda^N \int \lbrk{ \varphi(\lambda u_{\Omega_1}(z)) -
      \varphi(     \lambda u_{\Omega_2}(z)) }^p d\,z \label{eq:zrfsC}\\
    &=\lambda^N  d_{p,\tilde \varphi}(\Omega_1, \Omega_2)^p
  \end{align}
  where to go from (\ref{eq:zrfsA}) to (\ref{eq:zrfsB}) we used the change of variable $x = \lambda z$ and the property of the distance function
  \begin{equation}
    u_{\lambda \Omega}(\lambda z) = \lambda u_{\Omega}(z)\label{eq:scala_u}
  \end{equation}
  to change (\ref{eq:zrfsB}) to (\ref{eq:zrfsC}).
  % \end{Remark}

\end{itemize}

\begin{Remark}\label{rmk:d1notcompact}
  The inequality (\ref{eq:d1bound}) easily implies that the balls of the
  distance $d$ in general are not compact sets.  Indeed it is enough
  to consider a set $\Omega$ and the following ball: 
  $\mathbb D = \{A ~|~ d(A,\Omega)\le  2r\}$
  with $r=\|v_{\Omega} \| _{L^p}$. Then the sequence: $\set{{\Omega + n \tau} }_{n
    \in \N}$ with $\tau \in \real^N \setminus \set{0}$ is contained in
  $\mathbb D$ and it does not have any convergent subsequence.
\end{Remark}

To continue with our study of $d$, we prove this fundamental inequality.
\begin{Lemma}[Local equiboundedness]
  \label{prp:loc_eq_com}
  % Fix $p\in [1,\infty]$, 
  There is a  continuous   and increasing function
  $b:\real^+\to\real^+$ with $b(0)=0$
  such that, for any $\Omega, \Omega'\in\calI$ satisfying 
  \[ \| v_\Omega - v_{\Omega'}\|_{L^p} < b(r), \]
  then $\Omega'\subset \Omega + D_r$. 
  \begin{proof}\em
    Set  $K\defeq \Omega + D_r$. 
    It is easy to check that 
    $$v_K(x) = \varphi\big(  ( u_\Omega(x) - r) ^+\big).$$

    To prove the proposition for $p\in [1,\infty)$, suppose that
    $x_0\in \Omega'$, but $x_0\not \in K$; for $y\in B(x_0,r/2)$
    recall the simple triangular inequality
    \[ u_{\Omega}(y) \ge r - |x_0-y| \ge |x_0-y| \ge u_{\Omega'}(y)\]
    hence
    \[ v_{\Omega }(y) \le \varphi( r - |x_0-y|) \le \varphi (|x_0-y|) \le
    v_{\Omega'}(y) \]
    
    \begin{eqnarray*}
      \| v_\Omega - v_{\Omega'}\|_{L^p} ^p &\ge& 
      \int_{B(x_0,r/2)} |  v_{\Omega'} - v_{\Omega} |^pdx  \ge \\ &\ge&
      \int_{B(x_0,r/2)} |\varphi (|x_0-y|) - \varphi( r - |x_0-y|)   |^pdx 
      = b(r)^p 
    \end{eqnarray*}
    where 
    \[      b(r)^p ~\defeq~ {\omega_N}N \int_0^{r/2}
    t^{N-1}( \varphi(t) - \varphi(r-t))^p dt\]
    where $\omega_N$ is the volume of the ball $B_1$.
    It is easy to prove that $b$ is  continuous and increasing 
    (by direct derivation); that $b(0)=0$ and that
    $\lim_{r\to\infty} b(r) = \| \varphi(|x|) \|_{L^p}$.

%     If  $\varphi(r)=\exp(-r)$ then \ammark{RIVEDERE}
%     \begin{eqnarray*}
%       b(r)^p
%       &=&    {\omega_N}N \int_0^{r/2}   t^{N-1}( e^{-t} - e^{t-r})^p dt
%       =    {\omega_N}N r^N \int_0^{1/2} s^{N-1}( e^{-rs} - e^{r(s-1)})^p ds\\
%       &=& {\omega_N}N r^N  ( 1 - e^{-r})^p \int_0^{1/2} s^{N-1} e^{- r s p} ds
%       \leq {\omega_N}N r^N  2^{1-N} ( 1 - e^{-r})^p \int_0^{1/2}  e^{- r s p} ds\\
%       &\leq& \frac{\omega_N N}{2^{N-1}  p  }  \cdot r^{N-1}.  
%     \end{eqnarray*}
%     Then 
%     \begin{equation}
%  b(r) \leq \rbrk{ \frac{\omega_N N }{2^{N-1}  p  } }^{1/p} \cdot r^{(N-1)/p}. 
%    \end{equation}
    The case $p=\infty$ is simpler:
    in this case we can note that 
    \[\|v_\Omega - v_{\Omega'}\|_\infty\ge 
    v_{\Omega'}(x_0)- v_{\Omega }(x_0) \ge \varphi (0)- \varphi( r)    \]
    and set $b(r)=\varphi (0)- \varphi( r)$.
  \end{proof}
\end{Lemma}

\begin{Corollary} \label{cor:bound_d1_dh}
  As a corollary we obtain that for $d(\Omega,\Omega')$ small enough,
  \[ d_H(\Omega,\Omega')\le  b^{-1}\Big( d(\Omega,\Omega') \Big) ~.\]
\end{Corollary}

\begin{Remark} \label{rmk:jalsds}
  The above does not hold for arbitrarily large distance
  $d(\Omega,\Omega')$: indeed,
  let   $\Omega=\{0\}$ and $\Omega_n=\{n e_1\}$: then
  $d(\Omega,\Omega_n)\to  2\|\varphi(|x|)\|_{L^p}$
  (as we mentioned in \eqref{eq:separ}). 
\end{Remark}

We can also obtain a converse inequality, as follows
\begin{Lemma} \label{prp:dh_d1}
  There is a family of continuous  functions $f_R:[0,1]\to\real^+$  
   with $f_R(0)=0$,  such that, for any $\Omega, \Omega'\in\calI$,
  if $\Omega$ has diameter $R$  and   $d_H(\Omega,\Omega')<1$, then
  \[ d(\Omega,\Omega')\le  f_R\Big( d_H(\Omega,\Omega') \Big) ~.\]
  \begin{proof}
    We provide the proof for $p<\infty$.
    Note that if $\Omega$ has diameter $R$  and   $d_H(\Omega,\Omega')<1$, then
    $\Omega'$ has diameter at most $R+2$. Up to translation, suppose that 
    $B_{2R+4}$ contains both $\Omega$ and $\Omega'$: then
    \[ v_\Omega(x),v_{\Omega'}(x) \le \varphi( (|x|-2R-4)^+ ) \]
    so 
    \[ \int_{\real^N\setminus B_r}  |v_\Omega(x)\vee v_{\Omega'}(x)|^p dx \le  a_R(r)\]
    where
    \[ a_R(r) \defeq
    \int_{\real^N\setminus B_r} \varphi( (|x|-2R-4)^+ )^p dx = {\omega_N}{N} 
    \int_r^\infty t^{N-1} \varphi( (t-2R-4)^+ )^p dt\]
    and note that $a_R(r)\to 0$  for $r\to\infty$.
    At the same time, let $l(r)=\sup_{[0,r]}|\varphi'|$:    then 
    \[ \forall x\in B_r, \quad |v_\Omega(x) -v_{\Omega'}(x)|
    \le  l(r+4+2R) |u_\Omega(x)-u_{\Omega'}(x)|\]
    so 
    \begin{eqnarray*}
    \int_{B_r}  |v_\Omega(x)-v_{\Omega'}(x)|^p dx 
    \le \omega_N r^N l(r+4+2R)^p \sup_{x\in B_r}|u_\Omega(x)-u_{\Omega'}(x)|^p
    \le\\ \le \omega_N  r^N l(r+4+2R)^p d_H(\Omega,\Omega')^p
  \end{eqnarray*}

    Summarizing,
    \begin{eqnarray*}
      d(\Omega,\Omega')^p =
      \int_{\real^N\setminus B_r}  |v_\Omega(x)-v_{\Omega'}(x)|^p dx +
      \int_{B_r}  |v_\Omega(x)-v_{\Omega'}(x)|^p dx \le \\ \le 
       a_R(r)+     \omega_N  r^N    l(r+4+2R)^p d_H(\Omega,\Omega')^p
    \end{eqnarray*}
    Let eventually
    \[ g_R(s) = \inf_{r\ge 2R+4}  [ a_R(r)+
    \omega_N  r^N  l(r+4+2R)^p ~s]\]
    and note that it is concave and that $\lim_{s\to 0} g_R(s)=0$;
    and let $f_R(s) =\sqrt[p]{g_R(s^p)}$.
  \end{proof}
\end{Lemma}

Combining the two lemmas  \ref{prp:dh_d1} and \ref{cor:bound_d1_dh},
we obtain that
\begin{Theorem}\label{prp:topol_equiv}
  The topology induced by $d$ over the space $\calI$ is equivalent to the one induced by $d_H$.
\end{Theorem}
This implies that all properties of the Hausdorff distance listed in
proposition \ref{prp:properties_in_dH} are valid for the distance $d$ as well.

\subsection{Completeness and compactness}
By prop.~\ref{prp:topol_equiv}, we know that
$(\calI,d)$ is locally compact.

We now prove that it is complete:
\begin{Proposition}[Completness] The space $(\calI, d)$ is complete.
  \label{prp:d_complete}
  \begin{proof}
    Let $\Omega_n$ be a Cauchy sequence; this means that, 
    %setting    $v_n\defeq v_{\Omega_n}$,
    $\set{v_{\Omega_n}}_n \subset N_c$ is a Cauchy sequence: since 
    $L^p$ is complete, $v_{\Omega_n}\to g$ in $L^p$. It is well known 
    (see e.g. thm IV.9 in \cite{Brezis}) that,
    up to subsequence that we indicate with $\set{v_k}_k$, there is
    also convergence   $v_k(x) \to g(x)$ for almost all
    $x$; let $u_k(x)\defeq \varphi^{-1} v_k(x)$
    and $u=\varphi^{-1} g$; then $u_k(x) \to u(x)$ on a dense
    subset, so by the lemma \ref{lem:conv_punt}, $u=u_\Omega$
    where $\Omega\defeq\set{ u =0}$.
  \end{proof}
\end{Proposition}

Summarizing, this and \ref{prp:topol_equiv}
imply that $N_c$ is a complete (that is, closed)
and locally compact subset of $L^p$.

\begin{Remark}
  The above implies an interesting property of the subset $N_c$
  of $L^p$: it admits a small neighbourhood $U$ on $L^p$ such that,
  for $f\in U$, there is at least a $v\in N_c$ providing
  the minimum of the distance $\inf_{v\in N_c}\|f-v\|$.
  As far as we know, this minimum may fail to be 
  unique.  \ammark{[indagare ???]} 
\end{Remark}

\subsection{Shape analysis}
The family of distances is suitable for Shape Analysis:
we can indeed prove
\begin{Proposition}
  Let $G = {\mathcal O}(N) \ltimes \real^N$ be the Euclidean group of
  rotation and translation; as in \eqref{eq:eq:inv_quot}, we can define the
  quotient metric by
 
 \begin{equation} \label{eq:quotdistdef}
  d_q( [A], [B] ) = \inf_{g\in G} d(g A,  B).
  \end{equation}
  Then the above infimum is a minimum; so
  $d_q ([A], [B])>0$ when $[A] \neq [B]$.
  \begin{proof}
    Choose     a minimizing sequence
    $\set{g_n = (R_n, T_n) }_{n \in \N}$, that is
    \[ \inf_{g\in G} d(g A,  B) = \lim_{n\to \infty} d( g_n A , B) = \lim_{n\to \infty} d( R_n A + T_n , B).\] 
    Then $\set{T_n}_{n \in \N}$ must be bounded;
    in fact, let us assume by contradiction that $|T_n| \to \infty$, 
    then by \eqref{eq:d1bound} we would have that
		\[ d( A , B) < \| v_{A} \|_{L^p} + \| v_{B} \|_{L^p}  \]
	       and by (\ref{eq:separ}) that
  	\[ \lim_{n \to \infty} d(R_n A + T_n,B) =  \| v_{A} \|_{L^p} + \| v_{B} \|_{L^p},\]
  	so $\set{g_n}$ is not a minimizing sequence. 
  	This contradiction is generated by the assumption that $\set{T_n}$ is unbounded;
  	then the translation part of every minimizing sequence of (\ref{eq:quotdistdef}) must be bounded.
  	By compactness we have that there exists a limit transformation $g = (R,T) \in K$ such that $g_n \to g$ and by continuity
  	of $d( f A,B)$ with respect of $f \in G$, we have that $d( g A,B) =
        d^g( [A],[B])$.
  \end{proof}
\end{Proposition}

\subsection{$d^g$ and geodesics}
In this section we restrict $p\in(1,\infty)$.

Unfortunately $d=d_{p,\varphi}$ is not path--metric:
\begin{Proposition}\label{prop:dg_neq_d}
  Given any two $A,B\in\calI$ with $A \neq B$
  \begin{itemize}
  \item  then there is at most one $\lambda\in (0,1)$
    such that $\lambda v_A + (1-\lambda) v_B \in N_c$
    \ammark{\footnote{\tt forse: for all $\lambda\in (0,1)$, 
        $\lambda v_A + (1-\lambda) v_B\not\in N_c$}}
  \item consequently, by thm.~\ref{thm:geodetiche_tese},
    we have that    $d(A,B)<d^g(A,B)$.
  \end{itemize}
  \begin{proof}
    It is immediate to show that $f_\lambda =  \lambda v_A + (1-\lambda) v_B$ 
    assumes the value $\varphi(0)$ only on the intersection of the two sets $A \cap B$
    for any $\lambda \in (0,1)$. 
    Then $f_\lambda \in N_c$ implies that $f_\lambda = v_{A \cap B}$.
    Let $x \in (A \cap B)^c$ such that $v_A(x) \not= v_B(x)$. 
    We have that $f_{\lambda_1}(x) \not= f_{\lambda_2}(x)$ 
    for any $\lambda_1, \lambda_2 \in (0,1)$ 
    and $\lambda_1 \neq \lambda_2$. 
    Then there is at most one $\lambda \in (0,1)$ 
    such that $f_\lambda(x) = v_{A \cap B}(x).$ 
  \end{proof}
\end{Proposition}

So, to prove that the metric $d$ admits minimal geodesics,
we have to study $d^g$ as well;
to this end, we prove two results.

\begin{Proposition} \label{prop:Lip_arc_conn}
  If
  \begin{equation}
    \varphi'(|x|)\in L^p(\real^N)  \label{eq:phi'_Lp}
  \end{equation}
  then the space $(\calI,d_{p,\varphi})$ is Lipschitz--arc
  connected.
  \begin{proof}
    Indeed, let $\gamma(t)= t \Omega$ be the path that rescales $\Omega$
    to the singleton $\{0\}$; we prove that $\gamma$ is Lipschitz.

    It is not difficult to prove that the map 
    $(t,x)\mapsto u_{t\Omega}(x)$ is jointly Lipschitz. Then
    $u_{t\Omega}(x)$ is differentiable at
    almost all $t,x$, and fix such a $t,x$; note that
    \[ u_{t \Omega}(x) = t u_{\Omega}\left(\frac x t\right)\] (as in
    eqn.~\eqref{eq:scala_u}); hence, taking derivatives w.r.t.  $x$ we
    obtain
    \[ \nabla u_{t \Omega}(x) = \nabla u_{\Omega}\left(\frac x
      t\right)\] while taking derivatives w.r.t.  $t$ we obtain
    \[ \partial_t u_{t \Omega}(x) = u_{\Omega}\left(\frac x t\right)
    -\frac 1 t \langle \nabla u_{\Omega}\left(\frac x t\right) \cdot x
    \rangle = \frac 1 t \left( u_{t\Omega} (x) - \langle \nabla
      u_{t\Omega}(x ) \cdot x \rangle \right)~.\] 

    Suppose now that $x\not\in t\Omega$
    and let  $y\in t\Omega$ be a minimum distance point from $x$: then
    \[ u_{t \Omega}(x) = |x-y| ~~,~~ \nabla u_{t \Omega}(x) =
    \frac{x-y}{|x-y|} \] so
    \begin{eqnarray}
     \partial_t u_{t \Omega}(x) = \frac 1 t \left( |x-y| -
       \left\langle \frac{x-y}{|x-y|} \cdot x \right\rangle \right)
     = \nonumber \\ = -\frac 1 {t |x-y|}\langle x-y \cdot y\rangle =-\langle
     \frac{x-y}{|x-y|} \cdot \frac y t \rangle\label{eq:partial_t_u}
   \end{eqnarray}
   so if $\Omega\subset
    B_r$ we obtain that $|\partial_t u_{t \Omega}(x) | \le r$.
    If instead $x\in t\Omega$  and
    $u_{t\Omega}(x)$ is differentiable at $x$ then  $\nabla u_{t\Omega}(x)=0$
    and $\partial_t u_{t \Omega}(x)=0$.

    To conclude (cf. \ref{thm:len_d_B}) we compute
    \begin{equation}
      \|\dot\gamma\|_{L^p}^p =\int |\varphi'(u_{t \Omega}(x))|^p
      |\partial_t u_{t \Omega}(x) |^p ~dx \le r^p \int |\varphi'(u_{t
        \Omega}(x))|^p dx    \label{eq:last_W1p}
    \end{equation}
    and we argument as in \ref{lemma:compact_in_Lp}.
    By Rem.~1.1.3 in \cite{AGS04:gradien_flows_in_metric},
    we conclude that $\gamma$ is Lipschitz.
  \end{proof}
\end{Proposition}

\begin{Remark}
  Asking that $\varphi$ satisfy both \eqref{eq:phi_Lp} and  \eqref{eq:phi'_Lp}
  is equivalent to asking that  $\varphi(|x|)\in W^{1,p}$.
\quad
  By using the equality in  \eqref{eq:last_W1p}
  and in \eqref{eq:partial_t_u}, it is possible to show
  that, for most compact sets, the rescaling is a  Lipschitz
  path if and only if  $\varphi(|x|)\in W^{1,p}$.
\end{Remark}

% \ammark{\scriptsize attenzione: H2 diceva
%  \[    \lbrk{ \varphi(r + \lambda) - \varphi(r) } \leq \lambda \varphi(r), \ \ \ \forall r, \lambda  \geq 0.\]
% cioÃ¨
% \[    \varphi(r) -  \varphi(r + \lambda)   \leq \lambda \varphi(r), \]
% che Ã¨ ovvio per $\lambda \ge 1$; altrimenti
% \[    \varphi(r)(1-\lambda) \le  \varphi(r + \lambda), \]
% derivo in $\lambda =0$ e ho
% \[    -\varphi(r) \le  \varphi'(r ), \]
% dunque se H2 vale allora
%   \[ \varphi'(|x|)\in L^p\]
% e per questo motivo ho cancellato la vecchia proposizione.}

When $\calI$ is Lipschitz-arcwise connected, 
the induced metric $d^g=(d_{p,\varphi})^g$ is a finite
metric.

We can prove an equiboundedness result for $d^g$
(that is stronger than \ref{prp:loc_eq_com})
  \begin{Proposition}%[local equicompactness]
    \label{prp:equi_dg}
    Fix a compact nonempty set $\Omega$, and $r>0$;
    then there is a $K$ compact large such that for any
    closed set $\Omega'$ satisfying $d^g(\Omega,\Omega')< r$,
    then $\Omega'\subset K$.

    \begin{proof}
      Let $b(r)$ be  defined in \ref{prp:loc_eq_com}.
      Let $d^g(\Omega,\Omega')< r$, and
      $\gamma:[0,1]\to N_c$ be a Lipschitz path (of constant $L$)
      connecting $\gamma(0)=\Omega$ to $\gamma(1)=\Omega'$
      such that 
      \[\len \gamma \le d^g(\Omega,\Omega')+1\]
      up to reparametrization, we also assume that  $L\le r+2$.
      Let $n$ be large,     so that  $(r+2)/n\le b(r)$,
      and let $K=\Omega + D_{rn}$
      (note that $n$ only depends on $r$).
      Let $A_i=\gamma(i/n)$ for $i=0,\ldots, n$; we know that
      \[ d( A_i,A_{i+1}) \le d^g( A_i,A_{i+1}) \le L/n < (r+2)/n\le b(r)\] 
      since $\gamma$ is L-Lipschitz; so we
      apply recursively the proposition
      \ref{prp:loc_eq_com} on each $A_i$:  we obtain that, 
      \[ A_{i+1}\subset A_i+ D_r \]
      hence    $\Omega'\subset \Omega + D_{rn}=K$. 
    \end{proof}
  \end{Proposition}

\ammark{\scriptsize Corollario: dato $n$ t.c. $r/n\le b(r)$
  allora $ d^g(A,B)\le r \Longrightarrow  d_H(A,B)\le r n $
  ma questo porta a poco perche'
  $ d^g(A_n,A)\to 0 \Longrightarrow d(A_n,A)\to 0 \Longrightarrow
  d_H(A_n,A)\to 0$}

The above results have many interesting consequences:
\begin{Theorem}\label{thm:exists!}
  if \[ \varphi'(|x|)\in L^p\]
  then for any $\rho>0$,
  \[\mathbb D^g(A,\rho)\defeq \{ A ~|~ d^g(A,B)\le \rho\}\] is compact
  in the $(\calI,d)$ topology; so
  \begin{itemize}
  \item we obtain by Prp.~\ref{prp:D_compatto} that minimal geodesics
    do exist;
  \item and by  \ref{prp:exist_dgba} that the 
    \emph{Geodesic Distance Based Averaging}
    \begin{equation}
      \bar{A}=\argmin_A \sum_{j=1}^n d^g(A,A_j)^2 \label{eq:gdba_shape}
    \end{equation}
    of any given  collection $A_1,\ldots A_n$ exists.
  \end{itemize}
\end{Theorem}

\subsection{Variational description of geodesics}
In this section we restrict $p\in(1,\infty)$.
If $\gamma(t)$ is a Lipschitz path in $N_c$, then it is associated to
a function  $f(t,x)=v_{\gamma(t)}(x)$.
%; we state those results
\begin{Proposition}
  Suppose that $t\mapsto f(t,\cdot)$ is a Lipschitz path from $t\in
  [0,1]$ to $L^p(\real^N)$; then, by \ref{Radon-Nikodym}, for almost
  all $t$, $f$ admits strong derivative $\frac{df}{dt}$ in $L^p(\real^N)$
  (as was defined in eqn.~\eqref{eq:w_der}). Moreover
  \begin{itemize}
  \item $f$ admits weak partial derivative $\partial_t f$,
    and  $\partial_t f=\frac{df}{dt}$ for almost all $t$.
  \item If $f$ admits a pointwise partial derivative $h$ for almost all
    $t,x$, then \hbox{$\partial_t f=h$}.
  \end{itemize}
  \begin{proof}
    We extend %$\gamma(t)=\gamma(1)$ and
    $f(t,x)=f(1,x)$ for $t>1$, and %$\gamma(t)=\gamma(0)$ and
    $f(t,x)=f(0,x)$ for $t<0$; note that the extended $f(t,\cdot)$ is
    still Lipschitz in $L^p(\real^N)$; then we define
    \[g_\tau(t,x)\defeq\frac{f(t+\tau,x)-f(t,x)}\tau\]
    so 
    \[\|g_\tau(t,x)\|_{L^p(\real^N)}\le c\]
    where $c$ is the Lipschitz constant of $f(t,\cdot)$; hence
    \[\int_0^1 \int_{\real^N}|g_\tau(t,x)|^p~ dxdt\le c^p    \]
    This means that the family $g_\tau$ is bounded in 
    $L^p([0,1]\times\real^N)$, so we can find a sequence $\tau_n\to 0$ 
    such that $g_{\tau_n}\to w$ weakly, i.e.
    \[\lim_n\int_0^1\int_{\real^N} g_{\tau_n}(t,x) \psi(t,x)~dxdt=
    \int_0^1\int_{\real^N} w(t,x) \psi(t,x)~dxdt\]
    for all $\psi\in C^\infty_c([0,1]\times\real^N)$.  But 
    \[\int_0^1\int_{\real^N} g_\tau(t,x) \psi(t,x)~dxdt=
    \int_0^1\int_{\real^N} f(t,x) \frac{\psi(t-\tau,x)-\psi(t,x)}\tau ~dxdt \]
    hence
    \[\lim_n\int_0^1\int_{\real^N} g_{\tau_n}(t,x) \psi(t,x)~dxdt=
    -\int_0^1\int_{\real^N} f(t,x) \partial_t\psi(t,x)~dxdt\]
    by dominated convergence, 
    so we conclude that $f$ admits weak derivative, and the derivative
    is $w$. The relationship
    \eqref{eq:gamma_a_b} in  ${L^p(\real^N)}$, that is
    \[f(b,\cdot )  -f(a,\cdot ) = \int_a^b  \frac{df}{dt}dt\]
    implies that 
    \[\int_a^b \xi \frac{df}{dt}~dt= -\int_a^b   \frac{d\xi}{dt} f ~dt\]
    for all $\xi\in C^\infty_c([0,1])$;    but then setting
    $\psi(t,x)=\xi(t)$, we obtain that $\frac{d f}{dt}=\partial_t f$.
  \end{proof}
  \ammark{rinforzare con il teorema di derivabilita di Sobolev}
\end{Proposition}
This means that, for almost all $t$,
we can represent the ``abstract'' derivative
$\frac{d\gamma}{dt}$ by means of the weak derivative
 $\partial_t f(t,\cdot)\in L^p(\real^N)$.

We  use this result and eqn.~\eqref{eq:len_b_dot_gamma}
to express the length:
\begin{equation}
  \len^d\gamma = \int_0^1 \|\dot \gamma(t)\| dt
  = \int_0^1 \| \partial_t f(t,x) \|_{L^p} dt \label{eq:len_der_t_f}
\end{equation}
So to find the minimal geodesic between two compact sets $A,B$,
we need to minimize the above, with
the constraint that $f(0,\cdot)=v_A$, $f(1,\cdot)=v_B$, 
and, for any fixed t, $\varphi^{-1}f(t,\cdot)$ is a 
distance function.

It is possible to prove 
(using a reparametrization lemma and H\"oelder inequality)
that the geodesic is also the minimum of the \emph{action}
\[ J(\gamma) 
= \int_0^1  \| \partial_t f(t,x) \|_{L^p}^p dt
= \int_0^1 \int_{\real^N} | \partial_t f(t,x) |^p dx dt  \]

Equivalently, setting $g(t,x)=u_{\gamma(t)}(x)$, to 
find geodesics we can minimize
\[ J(\gamma) 
= \int_a^b \int_{\real^N} | \varphi'(\varphi(g))\partial_t g(t,x) |^p dx dt  \]
with the constraint that $g(0,\cdot)=u_A$, $g(1,\cdot)=u_B$, 
and, for any fixed t, $g(t,\cdot)$ is a distance function.

\subsection{Tangent bundle}
Let $p\in (1,\infty)$.
We identify  $\calI$ with $N_c\subset L^p$,
as by remark~\ref{iso}.

Given a $v\in N_c$, let $T_v N_c\subset L^p$ be the contingent cone
\[ T_v N_c \defeq 
\set{ \lim_n t_n ( v_n-v ) ~|~ t_n>0, v_n\in N_c, v_n\to v }
=
\set{ \lambda \lim_n  \frac{ v_n-v }{\|v_n-v\|}_{L^p} ~|~ \lambda\ge 0, v_n\to v }
~~, \] 
where it is intended that the above limits are in the sense of
strong convergence in $L^p$.

According to theorem~\ref{thm:len_d_B} 
if $\gamma:[a,b]\to N_c$ is a Lipschitz curve
then $\dot \gamma(t) \in T_{\gamma} N_c $ for almost all $t$.

In the following example we write explicitly the element 
of the contingent cone relative to a particular curve.
\begin{Example}
  We fix $\Omega \in \calI$, and define 
  the fattening  $\Omega_t = \Omega + D_t$ for $t \geq 0$. 
  We are interested 
  in evaluating the derivative $\dot{\gamma}(t)$. As
  previously done, we use the fact that
  \begin{equation}
    u_{\Omega_t}(x) = (u_{\Omega}(x) - t)^+
  \end{equation}
  and note that this map is jointly Lipschitz in $(t,x)$:
  hence both $u_{\Omega_t}(x)$ are $v_{\Omega_t}(x)$ are 
  almost everywhere differentiable.
  The pointwise  derivative is given by:
  \begin{equation}
    w = \lim_{\tau \to 0} \frac 1 \tau \qbrk{ v_{\Omega_{t+\tau}} - v_{\Omega_t}} = 
    \begin{cases}
      -\varphi'(u_{\Omega}(x) - t)     & \mbox{for \ } x \not\in  \Omega_t, \\
      0                & \mbox{for \ } x \in  \interno\Omega_t. 
    \end{cases}	
  \end{equation}
  (note that the derivative may not exists for $x\in\partial\Omega_t$).
%   but, for $t>0$, $\partial\Omega_t$ 
%   is a $N-1$ Lipschitz submanifold, so it is negligible).
  %se e' convesso,  class $C^{1,1}$
  If $\varphi'(|x|)\in L^p$ then $w\in L^p$,
  and it can be shown  that 
  \[w = \lim_{\tau \to 0}
  \frac 1 \tau \qbrk{ v_{\Omega_{t+\tau}} - v_{\Omega_t}}\]
  in the ${L^p}$ sense; then $w$ is in the  contingent cone.
  In particular, by Rem.~1.1.3 in \cite{AGS04:gradien_flows_in_metric},
  we obtain that the curve $\gamma$ is Lipschitz for $t\in [0,T]$.
\end{Example}

Unfortunately the contingent cone is not capable of
expressing some shape motions
\begin{Example}\label{exa:removing_in_Lp}
  We consider  the \emph{removing} motion;
  %that we saw in  \ref{prp:properties_in_dH}.(\ref{item:removing_in_dH}).
  to simplify the matter, let $A$ be compact, and suppose that
  the origin $0$ is in the internal part of $A$; let 
  $ A_t\defeq A\setminus B_t$ be the removal of a small
  ball from $A$: then we can explicitly compute (for $r>0,s>0$ small)
  \begin{eqnarray*}
    \| v_{A_{r+s}}- v_{A_s} \|^p_{L^p}
    &=& {\omega_N}N \int_s^{r+s} t^{N-1}\Big(\varphi(0)-\varphi(s+r-t)\Big)^pdt
    +\\ &+&{\omega_N}N \int_0^s t^{N-1}\Big(\varphi(s-t) -\varphi(s+r-t)\Big)^pdt
    \le \\ &\le&     {\omega_N}  r^p  L^p  (r+s)^N
  \end{eqnarray*}
  where $L$ is the  Lipschitz constant of $\varphi(t)$ for small $t$,
  and we see that this motion is Lipschitz.
  If we try to compute 
  \[    \frac{ v_{A_t}-v_A }{\|v_{A_t}-v_A\|}_{L^p} \]
  we notice that  $v_{A_t}-v_A=0$ outside of $B_t$: so the limit 
  would be zero for $x\neq 0$. 
%   for the particular case $\varphi(t)=\exp(-t), N=2, p=2$
%   \[   \frac{\pi}{2}  \left(e^{-r} (4 s-6)-2 e^{-r-2s}+e^{-2s}+ 
%     2r^2+(4 s-6) r-4 s+6+e^{-2 (r+s)}\right)\]
%   so that the metric derivative is
%   \[  \frac{1}{2} e^{-2 (r+s)} \left(e^{r+2 s} (3-2 s)+e^r+e^{2 (r+s)} (2 r+2
%     s-3)-1\right)\]
\end{Example}

\subsection{Riemannian metric}
\label{sec_riemannian}
Let now $p=2$.  
The set $N_c$ may fail  to be a smooth submanifold of $L^2$;
yet we will, as much as possible, pretend that it is, in order to
induce a sort of ``Riemannian metric'' on $N_c$ from the standard $L^2$ metric.

We  define the ``Riemannian metric'' on $N_c$ simply by
\[ \langle h,k \rangle \defeq \langle h,k \rangle_{L^2} \]
for $h,k\in T_v N_c$
and correspondingly a norm by 
\[|h| \defeq \sqrt{\langle h,h \rangle}\]
\begin{Proposition}
  We will also argue that the distance induced by this
  ``Riemannian metric'' coincides with the 
  geodesically induce distance $d^g$. Indeed
  %as in Thm.~\ref{thm:len_d_B}
  let $\gamma:[a,b]\to M$ be a Lipschitz curve  in $N_c$;  we may
  define the ``Riemannian length'' of the curve
  \[\len^R \gamma \defeq \int |\dot\gamma| ds\]
  Then we define the ``Riemannian distance'' $d^R(x,y)$ as the infimum
  of $\len^R\gamma$ for all $\gamma$ connecting $x$ to $y$.
  But by eqn.~\eqref{eq:len_b_dot_gamma}, $\len^R\gamma=\len\gamma$
  and $d^R = d^g$.
\end{Proposition}

\subsection{Example: smooth  convex sets}
\def\derpar#1{\partial_{#1}}

We propose, as an example, an explicit computation
of the Riemannian Metric.
We fix $p=2$,
% $\varphi(t)=\exp(-t)$,
$N=2$.
\quad
Let $\Omega\subset \real^2$ be a convex set with smooth boundary;
let $y(\theta):[0,L]\to \partial\Omega$ be a parametrization
of the boundary, $\nu(\theta)$  the unit vector normal  to $\partial\Omega$
and pointing external to $\Omega$: then the following
\emph{``polar''} change of coordinates holds:
$$\psi: \real^+ \times [0,L] \to \real\setminus \Omega \quad,\quad
\psi(\rho, \theta) =y(\theta) + \rho \nu(\theta)$$

We suppose that $y(\theta)$ moves on $\partial \Omega$
in anticlockwise direction; so 
\[ \nu= J \derpar s y \quad,\quad \derpar{ss} y = -\kappa \nu 
\quad,\quad  \derpar{s} \nu = \kappa \derpar s y \]
where $J$ is the rotation matrix
(of angle $-\pi/2$),  $\kappa$ is the curvature,
 and $\derpar s y$ is the tangent vector
(obtained by deriving $y$ with respect to arc parameter).

We can then express a generic integral through this change of coordinates as
$$\int_{\real^2\setminus \Omega} f(x) ~dx = 
\int_{\real^+} \int_{\partial{\Omega}} 
f( \psi(\rho, s) ) |1 + \rho \kappa(s)| ~d \rho d s$$
where $s$ is arc parameter, and $ds$ is integration in arc parameter.

\smallskip

We want to study a smooth deformation of $\Omega$, that we call  $\Omega_t$;
then the border $y(\theta,t)$ depends on a time
parameter $t$. %; we call $\alpha=\derpar t y$.
Suppose also that $\kappa(\theta)>0$, that is, that the set is strictly
convex: then for small smooth deformations, the set $\Omega_t$ will still be
strictly convex.
By deriving 
\[ \derpar t \derpar s y = \derpar s (\derpar t y) -
 \derpar s y \langle \partial_s y ,\partial_s (\derpar t y)\rangle =
 \pi_\nu( \partial_s (\derpar t y))\]
where 
\[\pi_\nu (w) \defeq w - \nu\langle \nu,w\rangle\]
is the projection of $w$ parallel to $\nu$. Supposing now that 
 $\rho=\rho(t)$ as well, we can express the point $\psi ( \rho,y)$
in a first order approximation as
\begin{eqnarray*}
  d \psi =
\Big( (\derpar t y) + \rho'  \nu + \rho  J\pi_\nu(\derpar s (\derpar t y)) \Big) dt+
 \Big( \derpar \theta y + \rho \derpar \theta \nu \Big) d\theta
\end{eqnarray*}
where moreover
\[  \Big( \derpar \theta y + \rho \derpar \theta \nu \Big) d\theta
=  \Big( \derpar s y + \rho \derpar s \nu \Big) ds =
\Big( 1+ \rho \kappa \Big)\derpar s y ds ~. \]

If $y(\theta,t),\rho(t)$ are expressing a constant point $x=\psi (
\rho,y)$, then $d\psi=0$; we apply scalar products w.r.t. $\nu$ and
$\derpar s y$ to the above relations
\[ \langle \nu, (\derpar t y)\rangle + \rho'=0\quad , \quad
 \langle \derpar s y, (\derpar t y)\rangle 
- \rho  \langle \nu , \derpar s (\derpar t y) \rangle dt+
  (1+ \rho \kappa  ) ds = 0~~.\]

Assuming that $(\derpar t y) \perp \derpar s y$, that is,
$(\derpar t y) = \alpha \nu $ with $\alpha=\alpha(t,\theta)\in\real$,
we obtain the relationships
\[ \rho' = -  \alpha \quad , \quad
\frac{ds}{dt} = \frac{\rho  \langle \nu , \derpar s (\alpha\nu) \rangle}
{(1+ \rho \kappa  )} =
\frac{\rho ~\derpar s \alpha }{(1+ \rho \kappa  )} ~. \]

\medskip

Now, for $x\not\in\Omega_t$, $u_{\Omega_t}(x) = \rho(t)$ hence
\[h_\alpha \defeq \derpar t v_{\Omega_t}(x)= -  \varphi'(u_{\Omega_t}(x)) \alpha \]
whereas $h_\alpha(x)=0$ for $x\in\interno \Omega_t$; so $h_\alpha$
is the vector in $T_vN_c$ that is associated to~$\alpha$.

Let us then fix two orthogonal smooth vector fields
$\alpha(s)\nu(s),\beta(s)\nu(s)$,
that represent two possible deformations of
$\partial{\Omega}$; those correspond to two vectors
$h_\alpha, h_\beta \in T_vN_c$;
so the Riemannian Metric that we presented in Sec.~\ref{sec_riemannian}
can be pulled back on $\partial{\Omega}$, to provide the metric
\begin{eqnarray*}
\langle \alpha,\beta\rangle &\defeq& \int_{\real^2} h_\alpha(x)h_\beta(x) dx=
\int_{\real^2\setminus \Omega} h_\alpha(x)h_\beta(x) dx = 
\\ &=& \int_{\partial{\Omega}}  \left[ \int_{\real^+}  (\varphi'(\rho))^2
 (1 + \rho \kappa(s)) ~d \rho\right] \alpha(s)\beta(s) d s 
\end{eqnarray*}
that is,
\begin{equation}
\langle \alpha,\beta\rangle \defeq  \int_{\partial{\Omega}}  (a + b \kappa(s) )\alpha(s)\beta(s) d s\label{eq:alfabeta}
\end{equation}
with
\[a = \int_{\real^+}  (\varphi'(\rho))^2 ~d\rho \quad,\quad
b =  \int_{\real^+}   (\varphi'(\rho))^2 \rho~d\rho ~.\]

\paragraph{Smooth sets}
If $\Omega$ is smooth but not convex, then the above formula holds up
to the cutlocus. We define a function $R(s):
[0,L]\to\real^+$ that spans  the cutlocus, that is, 
\[\text{Cut}= \{ \psi(R(s),s), s\in[0,L]\}~~.\]
$\psi$ is a diffeomorphism between the sets
\[ \{ (\rho,s)~  s\in[0,L], 0<\rho<R(s) \}
\leftrightarrow \real^2\setminus (\Omega\cup\text{Cut} )\]
moreover $R(s)$ is  Lipschitz
(by results in \cite{JI-MT:Lip},\cite{li-nirenberg03:finsl_hamil_jacob}).

In this case the metric has the form
\[ \langle h,k\rangle = 
\int_{\partial{\Omega}}  \left[ \int_0^{T(s)}  (\varphi'(\rho))^2
 (1 + \rho \kappa(s)) ~d \rho\right] \alpha(s)\beta(s) d s \]

\section{Other Banach--like metrics of shapes}
The paradigm that we presented in the previous section may be exploited
in other similar ways;
to conclude the paper, we shortly present some different embeddings
(leaving to a future paper the detailed study of their properties).

\subsection{Signed distance based representation}
We may use the signed distance function 
$b_A$, that was defined in \eqref{eq:signdistfun},
to define a metric of shapes:
\[d'(A,B) \defeq \|\varphi(b_A) - \varphi(b_B)
\|_{L^p(\real^N)} \]
in this case, we require that the function
$\varphi:\real\to(0,\infty)$  is monotonically
decreasing and of class $C^1$, and such that
\begin{equation}
  \varphi(|x|-t)\in L^p(\real^N) ~~~ \forall t.
\end{equation}

The resulting metric is slightly stronger than the one we studied in the
preceding sections; in particular,
\begin{Remark}
  Let $\mathcal F$ be the class of all finite subsets of $\real^N$; this
  class is dense in $\calI$ when we use the metric $d_{p,\varphi}$, or
  the Hausdorff metric; but it is not dense when we use the metric
  $d'$.
\end{Remark}

\subsection{$W^{1,p}$ metrics}
Another interesting choice of metric is obtained by embedding 
the representation in $W^{1,p}$, for $p\in(1,\infty)$

We require  that
$\varphi:[0,\infty)\to(0,\infty)$ 
 be Lipschitz, $C^1$ and monotonically decreasing,
and $\varphi(|x|)\in W^{1,p}(\real^N)$; 
for the case $p<\infty$ we are equivalently asking that
\[\int_0^\infty t ^{N-1}( \varphi(t)^p +
|\varphi'(t)|^p) ~dt <\infty\]
and this implies that 
$\lim_{t\to\infty} \varphi(t)=0=\lim_{t\to\infty} \varphi'(t)$.

\ammark{vorrei    levare questa ipotesi....}
We add one last hypothesis,
we   assume that there is a $T>0$ s.t.
$\varphi(t)$ is convex for $t\in [T,\infty]$.
\begin{Proposition}
  For any $A$ compact we have $v_A\in  W^{1,p}(\real^N)$.
  \begin{proof}
    We already know by \ref{lemma:compact_in_Lp} that $v_A\in  L^{p}(\real^N)$.

    By hypotheses above, $v_A$ is Lipshitz; and then, for almost all $x$,
    $\nabla v_A= \varphi'(u_A)\nabla u_A$; where $|\nabla u_A|=1$
    for almost all $x\not\in A$, while $\nabla u_A=0$
    for almost all $x\in A$.
    We also know 
    that when $t>T$, $\varphi'(t)<0$, $\varphi'$ is increasing and
    $\varphi'(t)\uparrow 0$.
    
    Let $R>0$ be large so that $A\subset B_R$, then
    \[ u_A(x) \ge |x|-R\]
    and then when $|x|\ge R+T$ we obtain that 
    \[ \varphi'(u_A(x)) \ge \varphi' (|x|-R)\]
    that is
    \[ \int_{\real^N\setminus B_{R+T}}|\varphi'(u_A(x))|^p~dx \le 
    \int_{\real^N\setminus B_{R+T}} |\varphi' (|x|-R)|^p~dx    <\infty~~.\]
    At the same time, since $v_A$ is Lipschitz, 
    then $\int_{B_{R+T}} |\nabla v_A|dx$ is finite.
  \end{proof}
\end{Proposition}

\begin{Definition}
  Given $A,B\in\calI$, we define
  \[d_{1,p,\varphi}(A,B) \defeq \|\varphi(u_A) - \varphi(u_B)
  \|_{W^{1,p}(\real^N)} \]
\end{Definition}

We just state a simple property of this metric:
\begin{Proposition}
  Let again $\mathcal F$ be the class of all finite subsets of $\real^N$:
  this class is dense in $\calI$ if and only if $\varphi'(0)=0$.

  Indeed, fix $A$ compact; let $\set{ x_n}_n$ be a dense subset of $A$;
  let $A_k\defeq \set{x_k~|~k\le n}$ be a finite subfamily;
  if $\varphi'(0)=0$ then $A_k\to A$ according to $d_{1,p,\varphi}$.

  If $\varphi'(0)<0$  it is easy to find examples where this does not hold:
  let $N=1$, $A=[0,1]$, then $\int_0^1 |v_{A_k}'(t)|^pdt\to |\varphi'(0)|^p$.
\end{Proposition}

\section*{Conclusions}
We have studied a metric space of shapes $(\calI,d_{p,\varphi})$; 
this space has a ``weak distance'', in that it has many compact sets,
and geodesics do exist;
but it can be associated in some cases to a smooth Riemannian metric,
as we saw in eqn.~\eqref{eq:alfabeta}.
% so it has the potential to overcome the tradeoff that we
% discussed in \ref{sec:antinomy}
% \footnote{and this is potentially eve}. magari per W1p
Moreover, by the properties that we saw in sec.~\ref{sec:emb}
(and in particular, by the properties of
$L^p$ spaces for $p\in(1,\infty)$ that we proved in 
Thm.~\ref{thm:geodetiche_tese}) we can also hope that geodesics
can be studied in the O.D.E. sense (altough possibly
in a very weak sense).

\medskip

As we saw in the last chapter, 
the representation/embedding paradigm can be exploited in many different
fashions; we just conclude with one last remark.
\begin{Remark}
  The embedding of $\varphi\circ u_A$ in $W^{2,p}$  is not feasible:
  if $A$ is smooth but is not convex, the second derivative
  of $u_A$ along the cutlocus is expressed by a measure 
  (see 4.13 in \cite{CM-ACM:HamJacEqu02})
  and then  $\varphi\circ u_A\not\in W^{2,p}$.
\end{Remark}

\tableofcontents

\bibliographystyle{plain}
\bibliography{../../mennbib}%,DistanzaL2}
% \bibliography{DistanzaL2}
\end{document}